\documentclass{article}
\usepackage[utf8]{inputenc}
\usepackage{fullpage}
\usepackage{hyperref}
\usepackage{xcolor}
\usepackage{amsmath,amsfonts,amssymb}
\usepackage{graphicx}%
\graphicspath{{New_plots/}}
\usepackage{subcaption}
\usepackage{bbm}
\usepackage{array}
\usepackage{graphicx}%
\newtheorem{theorem}{Theorem}

\newtheorem{remark}[theorem]{Remark}
\usepackage{algorithm2e}
\usepackage{algpseudocode}

\newcommand{\rd}{\mathrm{d}}
\newcommand{\loc}{\text{loc}}
\newcommand{\equi}{\text{equi}}
\newcommand{\Kn}{\mathsf{Kn}}

\newcommand{\RR}{\mathrm{R}}
\newcommand{\LL}{\mathrm{L}}
\newcommand{\A}{\mathrm{A}}
\newcommand{\interior}{\mathrm{in}}
\newcommand{\Null}{\mathrm{Null}}

\begin{document}

\title{Second-order diffusion limit for the phonon transport equation -- asymptotics and numerics}
\author{Anjali Nair\thanks{University of Wisconsin-Madison, Mathematics department {nair25@wisc.edu}} and Qin Li\thanks{University of Wisconsin-Madison, Mathematics department and Wisconsin Institute for Discovery, {qinli@math.wisc.edu}} and Weiran Sun\thanks{Simon Fraser University, Mathematics department,weirans@sfu.ca}}
\maketitle

\begin{abstract}
We investigate the numerical implementation of the limiting equation for the phonon transport equation in the small Knudsen number regime. The main contribution is that we derive the limiting equation that achieves the second order convergence, and provide a numerical recipe for computing the Robin coefficients. These coefficients are obtained by solving an auxiliary half-space equation. Numerically the half-space equation is solved by a spectral method that relies on the even-odd decomposition to eliminate corner-point singularity. Numerical evidences will be presented to justify the second order asymptotic convergence rate.
\end{abstract}

\section{Introduction}
Heat is a physical phenomenon that describes temperature fluctuation. The classical description for heat conductance is the simple heat equation, a parabolic equation when time is present, or an elliptic type in the steady state. The derivation of the heat equation is based on Fourier law that states the heat flux proportionally depends on the temperature gradient, so the bigger temperature fluctuation leads to stronger heat flux. This law is an observational fact but is not derived from the first principle.

In modern physics, derived from the first principle, it was discovered that this Fourier law may not be accurate. The underlying physics model for heat propagation should be characterized by the phonon transport equation. It describes the dynamics of phonons, the microscopic quanta that propagates heat energy. This phonon transport equation, builds upon the Wigner transform from quantum mechanics, and hence first principle, is a mesoscopic description that follows the statistical mechanics derivation. Based on this phonon transport equation, at the correct scaling, the heat equation is then rediscovered as the associated macroscopic limiting system when the system is ``large" enough so that the mesoscopic fluctuation can be ignored. This partially justifies the validity of the heat equation that has been traditionally used as a model equation~\cite{peraud2016extending}. Some experimental and numerical results are available in~\cite{hu2015spectral, minnich2011quasiballistic,  hua2014transport}.

In this article, we study this asymptotic relation between the phonon transport equation and the diffusion heat equation, and we pay special attention to the boundary effect. When no physical boundary is present, the derivation from one equation to the other is a rather straightforward process. When physical boundary is present, however, boundary layers emerge adjacent to the boundaries. They are used to damp the fluctuations at the mesoscopic level close to the physical boundaries that are inconsistent with the limiting equation. At the macroscopic limit, this means the boundary condition for the limiting diffusion equation needs to be fine-tuned to reflect such perturbation.

The study of kinetic layers and their effects has attracted quite some attention during the past decade~\cite{ degond2005smooth, golse2003domain, lemou2012micro, crouseilles2003hybrid}. One example that we follow in this paper is~\cite{li2015diffusion}, in which, the authors studied the neutron transport equation and derived the Dirichlet boundary condition, as the leading order approximation to the boundary condition for the limiting equation. The idea is to perform asymptotic expansion, and separate the studies of the interior and boundary layers. The characterization of boundary layers will be rephrased into a half-space problem, with its information at the infinite point translated to the boundary condition for the interior part of the domain. The well-posedness of the half-space equation, and an efficient numerical solver that achieves the spectral accuracy are discussed in~\cite{li2017convergent}. Earlier attempts can be found in~\cite{klar96, BM:02, CoronGolseSulem:88,CD:13}.

In this article, we follow the steps conducted in~\cite{golse2003domain,li2015diffusion,peraud2016extending} and study the boundary layer effect of the phonon transport equation with either the incoming or the reflective boundary conditions, with a goal of recovering a higher order asymptotic approximation. The associated half-space equation that characterizes the boundary layer evolution will be derived, and a numerical scheme based on spectral method will be designed to determine the Robin boundary value for the limiting equation. The whole scheme provides a second order asymptotic approximation to the original phonon transport equation. The current paper mainly differs from~\cite{li2015diffusion, peraud2016extending} in two aspects: 1. In the earlier papers, the authors derived the Dirichlet boundary condition, as the limiting incoming boundary condition. This only captures the leading order information asymptotically. In this paper, we extend the asymptotic expansion to the next order, and for a larger class of boundary conditions (both incoming and reflective), we derive a Robin type boundary condition for the macroscopic limit, increasing the asymptotic convergence by another order. This change mostly brings some technical difficulties, but the structure of derivation is unchanged; 2. In~\cite{CD:13,li2015diffusion} where the authors study the neutron transport equation, the frequency domain is not considered, based on the assumption that particles at different energy levels do not exchange information. We remove this assumption in this paper, and allow particles at all energy levels to interact, so the distribution converges to an equilibrium in both velocity and frequency domains. One immediate difficulty it brings is that the Knudsen number differ at different energy levels, meaning the equilibrium is achieved at different ``time." A universal Knudsen number is thus defined to unify the convergence rate. The process for deriving both the limiting equation and the boundary condition is accordingly altered to accommodate this change.

This paper is arranged in the following way. In Section~\ref{sec:eqn_asymptotics}, we discuss the phonon transport equation and derive its asymptotic limit. In Sections~\ref{sec:Boundary layer} and~\ref{sec:Diffusion limit}, we formulate the boundary layer correction and use this to compute the diffusion equation. Finally, we present the numerical examples in Section~\ref{sec:Numerical results}.

\section{Equation and its asymptotic limit}\label{sec:eqn_asymptotics}
In this section we study the equation and derive its asymptotic limit under the assumption that there is no boundary, or the boundary condition is compatible, to remove the existence of boundary layer effect. In particular, the equation is described in section~\ref{sec:equation}, and its second order asymptotic expansion is given in section~\ref{sec:asymptotics}

\subsection{Phonon transport equation}\label{sec:equation}
Phonon transport equation is a PDE model for describing the dynamics of phonon and heat propagation. It is a mescoscopic description and sits in the classical statistical mechanics framework. The equation writes as

\begin{equation*}
    \boldsymbol{v}_g\cdot\nabla_\mathbf{x}F=\frac{F^{\loc}-F}{\tau}\,,
\end{equation*}
where $F(x,v,\omega)$ is the phonon occupation number in the phase space $(\mathbf{x},\mathbf{v},\omega)\in\mathbb{R}^3\times \mathbb{S}^2\times\mathbb{R}^+$. Here $\mathbf{x}$ is the spatial variable, and $\boldsymbol{v}_g=\|\boldsymbol{v}_g\|\mathbf{v}$ is the phonon group velocity, where $\mathbf{v}\in\mathbb{S}^2$ presents the direction of the velocity, and $\omega$ is the photon frequency. Typically the group velocity $\|\boldsymbol{v}_g\|$ has the frequency dependence. $\tau$ is the phonon relaxation time. It also depends on $\omega$ and describes the amount of time for phonon at $\omega$ frequency to stablize to the equilibrium. $F^{\loc}$ defines the local equilibrium distribution. It depends on the local temperature $T$ via a Bose-Einstein distribution~\cite{hao2009frequency}:
\begin{equation}\label{eqn:equilibrium}
F^\loc = F^\equi_{T} = \frac{1}{e^{\hbar\omega/k_{\text{B}}T}-1}\,,
\end{equation}
where $T$, the local temperature is set to ensure the energy conservation, meaning:
\begin{align}
    \int\limits_{\mathbf{v},\omega}\frac{D\hbar\omega}{\tau}(F^{\loc}-F) \rd\omega\rd \mathbf{v}=0\,.
\end{align}
Here $\hbar$ is the rescaled Planck constant and $\hbar\omega$ is the amount of energy contained in one phonon at frequency $\omega$ and $D(\omega)$ is the phonon density of states. We note that due to the complicated dependence of $F^\loc$ on $F$, the equation is naturally nonlinear.

Suppose heat injected into the environment is not significant enough to deter the temperature drastically. Then one can linearize the system around the room temperature. To do so, we call $T^\ast$ the room temperature, and define the equilibrium at the room temperature:
\[
F^\ast = F^\equi_{T^\ast} = \frac{1}{e^{\hbar\omega/k_{\text{B}}T^\ast}-1}\,.
\]
Denote $F^d=F-F^\ast$. Considering $F^\ast$ has no $x$ dependence, we rewrite the original equation for $F^d$:
\begin{align}
    \boldsymbol{v}_g\cdot\nabla_\mathbf{x}F^d=\frac{F^{\loc}-F^d-F^\ast}{\tau}\,.
\end{align}
In this small fluctuation regime, the temperature change is small so that $F^\loc\sim F^\ast$, and we have the explicit formulation for temperature increase:
\begin{align}
F^\loc-F^d-F^\ast\approx\Delta T\frac{\partial F^{\loc}}{\partial T}\Big|_{T^\ast}-F^d\,,
\end{align}
where $\Delta T$ can be determined using the energy conservation. Setting $\int_{\mathbf{v},\omega}\frac{D\hbar\omega}{\tau}\left(\Delta T\frac{\partial F^\loc}{\partial T}|_{T^\ast}-F^d\right)\rd\mathbf{v}\rd\omega=0$, we have:
\begin{equation}\label{eqn:C_omega}
\Delta T=\frac{1}{C_\tau}\int\limits_{\mathbf{v},\omega}\frac{D\hbar\omega}{\tau}F^d \rd\omega\rd \mathbf{v}\,,\quad\text{with}\quad C_\tau=\int\limits_{\mathbf{v},\omega}\frac{C_\omega}{\tau}\rd\omega\rd \mathbf{v}\quad\text{and}\quad C_\omega=D\hbar\omega\frac{\partial F^{\loc}}{\partial T}\Big|_{T^\ast}\,.
\end{equation}

Assume plane-symmetry, the problem becomes pseudo-1D. For easier rescaling, we denote $f=F^d/(\frac{\partial F^{\loc}}{\partial T}\Big|_{T^\ast})$,  $v = v_1$ and rescale the spatial coordinate $x=\mathbf{x}_1/L$ by the characteristic length $L$, then the linearized equation becomes:
\begin{align}\label{eqn:f}
   v\partial_xf=\frac{\Delta T-f}{\Kn} = \frac{\mathcal{L}f-\alpha f}{\mathsf{Kn}} \,,\quad\text{with}\quad \alpha = 1+\alpha_0\langle\Kn\rangle^2\,,
\end{align}
where the Knudsen number is defined by the ratio of $\|\boldsymbol{v}_g\|\tau$ so that $\mathsf{Kn}(\omega)=\frac{\|\boldsymbol{v}_g\|\tau}{L}$. The physical meaning of the Knudsen number is still the distance for a phonon to travel before equilibriumized. The collision operator $\mathcal{L}$ is:
\begin{align}
    \mathcal{L}\phi=\langle\phi\rangle:=\frac{1}{C_\tau}\int\limits_{v,\omega}C_\omega/\tau\phi\rd v\rd\omega\,,
\end{align}
where we use the bracket notation $\langle\cdot\rangle$ to simplify the notation and the measure $\rd v\rd\rd\omega$ is normalized such that $\mathcal{L} (1) = 1$. We note that the Knudsen number depends on $\omega$, reflecting the fact that heat carried by phonons at different frequencies have various relaxation time and group velocity speed. To unify the derivation, we also define the averaged Knudsen number
\begin{align}
    \langle\Kn\rangle:=\frac{1}{C_\tau}\int\limits_{v,\omega}C_\omega/\tau\Kn\rd v\rd\omega\,.
\end{align}
Typically it is assumed that $\boldsymbol{v}_g$ and $\tau$ do not fluctuate severely with respect to the frequency, meaning that $\Kn$ at all $\omega$ have the same order of magnitude and thus $\langle\Kn\rangle\sim\Kn$. $\alpha$ is the total damping coefficient with $\alpha_0\geq 0$ being the extra attenuation. The order is chosen so to make the asymptotic limit unchanged.

\subsection{Classical asymptotic derivation}\label{sec:asymptotics}
When the domain $L$ is significantly larger than the average travel distance $\|\boldsymbol{v}_g\|\tau$, the Knudsen number $\Kn$ is quite small. In this case, one can derive the asymptotic limit of the phonon transport equation that characterizes its macroscopic behavior. Under the assumption that $\Kn$ and $\langle\Kn\rangle$ are at the same order of magnitude, we have the following theorem. 

\begin{theorem}
Let $f$ satisfy the Cauchy problem with the equation~\eqref{eqn:f}, and assume $\Kn\to 0$ (or $L\to\infty$), then in this limiting regime, $f(x,v,\omega)\to\rho(x)$ that satisfies:
\begin{equation}\label{eqn:rho_limit}
\frac{1}{3}\frac{\langle\Kn^2\rangle}{\langle\Kn\rangle^2}\partial_{xx}\rho +\langle\alpha_0\rangle\rho= 0\,.
\end{equation}
Moreover, one has $f(x,v,\omega)=\rho(x)-v\Kn\partial_x\rho +O(\Kn^2)$.
\end{theorem}
The theorem is a natural extension of the celebrated results in~\cite{BardosSantosSentis:84}. We opt not to repeat the proof but to give a formal derivation below.

Consider the following expansion in terms of $\langle \mathsf{Kn}\rangle$:
\begin{align}\label{eqn:f_expansion}
    f&=f_0+\langle \mathsf{\mathsf{Kn}}\rangle f_1+\langle \mathsf{Kn}\rangle^2f_2+f_r\,,
\end{align}
where $f_r$ stands for the remainder term. Inserting this expansion back in~\eqref{eqn:f} and comparing terms at different orders, we have:
\begin{itemize}
    \item[--]{At $\mathcal{O}\big(\frac{1}{\langle \mathsf{Kn}\rangle}\big)$:}
    \begin{align}\label{eqn:f_0_order}
    (\mathcal{L}-\mathbb{I})f_0&=0\,;
\end{align}
\item[--]{At $\mathcal{O}(1)$:}
\begin{align}\label{eqn:f_1_order}
    v\partial_xf_0&=\frac{\langle \mathsf{Kn}\rangle}{\mathsf{Kn}}(\mathcal{L}f_1-f_1)\,;
\end{align}
\item[--]{At $\mathcal{O}(\langle \mathsf{Kn}\rangle)$:}
\begin{align}\label{eqn:f_2_order}
        v\partial_xf_1&=\frac{\langle \mathsf{Kn}\rangle}{\mathsf{Kn}}(\mathcal{L}f_2-f_2)-\alpha_0 f_0\,.
\end{align}
\end{itemize}
Considering that $\mathcal{L}$ is an operator for $\omega$ and $v$, then noting~\eqref{eqn:f_0_order}, and that the null space of the operator $\mathcal{L}-\mathbb{I}$ is given by
\begin{equation}\label{eqn:null_L}
    \Null({\mathcal{L}-\mathbb{I}}) = \text{Span}\{1\}\,,
\end{equation}
we have $f_0$ independent of $\omega$ and $v$, so we set: $f_0(x,v,\omega)=\rho(x)$. Inserting this into~\eqref{eqn:f_1_order}, we have the solution to $f_1$ as $   f_1=-\frac{\mathsf{Kn}}{\langle \mathsf{Kn}\rangle}v\partial_x\rho$.
From the compatibility condition for $f_2$, we have
\begin{align*}
\frac{\langle v^2\Kn^2\rangle}{\langle\Kn\rangle^2}\partial_{xx}\rho +\langle\alpha_0\rangle\rho= 0\,.
\end{align*}
The constant $1/3$ comes from the integration of $v^2$ and it would take a different value in a different dimension. This derivation gives us the diffusion equation stated in the theorem. In the rigorous proof we need to give the bound to the remainder term $f_r$ and show its independence of $\Kn$. We omit the proof from here.

\section{Boundary layer correction}\label{sec:Boundary layer}

We now study the equation and its limit when physical boundaries are introduced. As discussed in the previous section, the asymptotic limit for the phonon transport equation can be derived when $f$ is approximately close to $\rho$, meaning asymptotically it loses its $(v,\omega)$ dependence, and becomes a constant for every different $x$.

When there is no extra boundary condition imposed on the physical boundaries, this is achieved by setting $\Kn$ extremely small, representing the long time large space limiting regime, as shown in the previous section. However, no matter how small $\Kn$ is, the boundary condition imposed on the physical boundary can significantly deter $f$ away from $\rho$. It takes a layer of the size of $\langle\Kn\rangle$ to damp the extra information in $f$ and reduce it to a constant $\rho$. As a consequence, the boundary information in $f$ is translated through the layer equation to become the boundary condition for $\rho$.



We study how to formulate boundary conditions on $\rho$ using the boundary conditions on $f$ in this section. To do so, we first reformulate the boundary layer equation into a half-space equation with some re-defined conditions. The two ends of the half-space equation are to connect the physical boundary and the interior condition.

To showcase our analysis in a compact manner, we restrict ourselves to the following equation:
\begin{equation}\label{eqn:f_original}
\begin{aligned}
\begin{cases}
    v\partial_x f&=\frac{1}{\Kn}\left(\mathcal{L}f-f\right)\\
    f(x=0,v,\omega)&=\phi,\quad v>0\\
    f(x=1,v,\omega)&=\eta(\omega)f(x=1,-v,\omega),\quad v<0
    \end{cases}
\end{aligned}
\end{equation}
meaning the phonon transport equation is supported in the domain $[0,1]$ and has an incoming left boundary and reflective right boundary. The reflection coefficient at $x=1$ has $\omega$ dependence. This equation not only serves as a toy problem for analyzing layer effect, but is also extracted from a practical experimental setup, see~\cite{hua2017experimental}, where the authors inject heat to two adjacent materials, and study the heat conductance at the interface of two solids.

As argued above, there will be thin layers in both left and right physical boundaries of the domain, and the layers are of width $\langle\Kn\rangle$. To proceed, we will write the approximated solution as
\begin{equation}\label{eqn:f_total}
f^\A = f^\LL + f^\interior + f^\RR\,,
\end{equation}
with the three terms taking care of the left boundary layer, interior solution and right boundary layer respectively. The superscript $\A$ stands for approximation. The hope is to derive proper equations for the three terms respectively so that $f^\A$ is asymptotically close to $f$. In~\cite{li2015diffusion} the equations were derived that obtained the leading order approximation. The goal is to modify the boundary conditions to achieve the higher order accuracy.

\begin{itemize}
    \item $f^\interior$: Following the procedure in~Section~\ref{sec:asymptotics}, it is straightforward to set
\begin{equation}\label{eqn:f_interior}
f^{\interior} = \rho-v \mathsf{Kn}\partial_x\rho\,,
\end{equation}
with $\rho$ satisfying~\eqref{eqn:rho_limit}. This is a 1D elliptic equation and one needs two boundary conditions to uniquely determine $\rho$.
\item $f^\LL$: Noting that the layer is of the size $\langle\Kn\rangle$, we first rescale the problem by defining $z=x/\langle \mathsf{Kn}\rangle$, then the rescaled equation becomes:
\begin{equation}
v\partial_zf^\LL=\frac{\langle \mathsf{Kn}\rangle}{\mathsf{Kn}}(\mathcal{L}f^\LL-f^\LL)\,.
\end{equation}
Furthermore, we need $f^\LL$ to have the eliminated effect in the interior, meaning we need $f^\LL(x) = 0$ for all finite $x\neq 0$. Calling the rescaling, we require:
\[
f^\LL(z\to\infty)=0\,.
\]
At the physical boundary $x=0$, we need the incoming boundary condition to be satisfied, meaning
\begin{equation}
f^\LL(z=0,v,\omega)+f^{\interior}(x=0,v,\omega)=\phi,\quad v>0\,.
\end{equation}
Calling~\eqref{eqn:f_interior}, this gives:
\begin{equation}
    \begin{aligned}
    f^\LL(z=0,v,\omega)&=\phi(v,\omega)-\rho(x=0)+v\mathsf{Kn}\partial_x\rho(x=0),\quad v>0\,.
    \end{aligned}
\end{equation}
We summarize everything to get the PDE $f^\LL$ satisfies:
\begin{equation}\label{eqn:f_L}
\begin{cases}
    v\partial_zf^\LL&=\frac{\langle \mathsf{Kn}\rangle}{\mathsf{Kn}}(\mathcal{L}f^\LL-f^\LL)\\
    f^\LL(z=0,v,\omega)&=\phi-\rho(x=0)+v\mathsf{Kn}\partial_x\rho(x=0),\quad v>0\\
    f^\LL(z\to\infty)&=0\,.
\end{cases}
\end{equation}
This is a half-space layer equation with Dirichlet boundary condition (or in-coming boundary condition). The condition itself is composed of $\phi$, the external condition, $\rho(x=0)$, a constant, and $v\mathsf{Kn}\partial_x\rho(x=0)$, a linear function in $v$.
\item $f^\RR$: The derivation for $f^\RR$ is the same. We consider the scaling $z=\frac{1-x}{\langle \mathsf{Kn}\rangle}, v\to-v$, then:
\begin{equation}\label{eqn:f_R}
\begin{cases}
v\partial_zf^\RR&=\frac{\langle \mathsf{Kn}\rangle}{\mathsf{Kn}}(\mathcal{L}f^\RR-f^\RR)\\
f^\RR(z=0,v,\omega)&=\eta f^\RR(z=0,-v,\omega)-(1-\eta)\rho(x=1)-(1+\eta)v\mathsf{Kn} \partial_x\rho(x=1),\quad v>0\\
f^\RR(z\to\infty,v,\omega)&=0\,.
\end{cases}
\end{equation}
This is a half-space layer equation with reflective boundary condition with a Dirichlet component. The reflective part takes a coefficient $\eta$ that has $\omega$ dependence, and the Dirichlet component is composed of $(1-\eta)\rho(x=1)$, a constant, and $(1-\eta)v\mathsf{Kn}\partial_x\rho(x=1)$, a linear function in $v$. Both components have $\omega$ dependence due to the involvement of $\eta$.

While the condition $f^\RR(z\to\infty,v,\omega)=0$ comes from the requirement that $f^\RR(x) = 0$ for all $x\neq 1$, the boundary condition for $f^\RR(z=0)$ is derived from the fact that $x=1$, according to~\eqref{eqn:f_original}:
\begin{align}
    f^{\interior}(x=1,v,\omega)+f^\RR(z=0,v,\omega)=\eta(\omega)[f^{\interior}(x=1,-v,\omega)+f^\RR(z=0,-v,\omega)]
\end{align}
and that $f^\interior$ is defined by~\eqref{eqn:f_interior}.
\end{itemize}

As a summary, the solution is presented by~\eqref{eqn:f_total} with the three terms solving~\eqref{eqn:f_interior},~\eqref{eqn:f_L} and~\eqref{eqn:f_R} respectively. The wellposedness of this whole system determines the value of $\rho(x=0)$ and $\rho(x=1)$.

We argue this decomposition gives a highly accurate approximation to the original equation set.
\begin{remark}
We do not provide the rigorous proof in this paper. The interested reader is referred to~\cite{golse_jin_levermore} for the parallel analysis. The statement is that, let $f$ solve~\eqref{eqn:f_original} and let $f^{\interior}$, $f^\LL$ and $f^\RR$ solve~\eqref{eqn:f_interior},~\eqref{eqn:f_L} and~\eqref{eqn:f_R} respectively. Then
\[
f^\A(x,v,\omega) = f^\interior(x,v,\omega) + f^\LL(\frac{x}{\epsilon},v,\omega) + f^\RR(\frac{1-x}{\epsilon},v,\omega)
\]
approximates $f$ with $O(\langle\Kn\rangle^2)$ accuracy:
\[
\|f-f^{\A}\|_{L^\infty(\rd{x}\rd{v}\rd\omega)}=O(\langle\Kn\rangle^2)\,.
\]
To show this, however one needs to asymptotically expand both the interior solution and the boundary equations to the higher orders, as was done in~\eqref{eqn:f_expansion}.
\end{remark}

We now examine the equation for $f^\LL$ and $f^\RR$ in more details. Indeed, to ensure the unique solution to~\eqref{eqn:f_L} and~\eqref{eqn:f_R}, one needs $f^\LL(z=0)$ and $f^\RR(z=0)$ to satisfy certain conditions. In particular, in~\cite{CoronGolseSulem:88}, the authors showed the wellposedness of a large class of linear kinetic equations. If applied here, this justifies that there is a unique solution to~\eqref{eqn:f_L} and~\eqref{eqn:f_R}. We cite the theorem here:

\begin{theorem}\label{thm:wellposedness}[Adaptation from Theorem 1.7.1 from~\cite{CoronGolseSulem:88} to fit the current setting.]
Let $\phi\in L^2(\rd{v}\rd\omega)$, then the solution to the following half-space equation:
\begin{equation}
    \begin{cases}
            v\partial_zf &=\frac{\langle\Kn\rangle}{\Kn} \left(\mathcal{L}f - f\right)\\
f(z=0,v>0,\omega) &= \phi
    \end{cases}
\end{equation}
has a unique solution if $f(z=\infty)\in\Null(\mathcal{L}-\mathbb{I})$.
\end{theorem}

We note that there is a very delicate difference between the conclusion of the theorem and equation~\eqref{eqn:f_L}-\eqref{eqn:f_R}. Take $f^\LL$ for example, when applied the theorem, we will have a unique solution for $f^\LL$ as long as $f^\LL(z\to\infty)\in\Null(\mathcal{L}-\mathbb{I})$. This is different from what we are looking for. As suggested in~\eqref{eqn:f_L}, we require $f^\LL$ to go to $0$ as $z\to\infty$. While it is true that $0\in\Null(\mathcal{L}-\mathbb{I})$, $f^\LL(z=0)$ has to satisfy certain conditions to make the two ends compatible. To do so, recall that the layer equation~\eqref{eqn:f_L} has three components in its incoming part, we accordingly set the solution to the summation of three pieces
\begin{equation}\label{eqn:f_L_decompose}
f^\LL=f_0-\rho(x=0)f_1+\partial_x\rho(x=0)f_2\,,
\end{equation}
with each piece $f_i, i=0,1,2$ solving the same equation but one component of the incoming data:
\begin{equation}\label{eqn:f_L_i}
\begin{cases}
v\partial_zf_i&=\frac{\langle \mathsf{Kn}\rangle}{\mathsf{Kn}}(\mathcal{L}f_i-f_i)\\
f_i(z=0,v,\omega)&=\psi_i\quad v>0\\
f_i(z\to\infty)&\in\Null(\mathcal{L}-\mathbb{I})
\end{cases}
\end{equation}
where $\psi_i$ takes the value of
\begin{equation*}\label{eqn:f0}
\psi_0(z=0,v,\omega)=\phi\,,\quad \psi_1(z=0,v,\omega)=1\,,\quad \psi_2(z=0,v,\omega)=v\mathsf{Kn}\,.
\end{equation*}
According to Theorem~\ref{thm:wellposedness}, we can find unique solutions for all $f_i$. Noting~\eqref{eqn:null_L}, we denote
\begin{equation}\label{eqn:b_i}
b_i=\lim_{z\to\infty}f_i\in\Null(\mathcal{L}-\mathbb{I})\,,\quad i = 0,1,2\,.
\end{equation}
Then naturally $f^\LL(z=\infty) = b_0 -b_1\rho(x=0) + b_2\partial_x\rho(x=0)$. Since we require, according to~\eqref{eqn:f_L_decompose}, that $f^\LL(z\to\infty)=0$, we have:
\[
b_0 -b_1\rho(x=0) + b_2\partial_x\rho(x=0) = 0\,.
\]
This becomes the requirement for $\rho$ at $x=0$. A Robin type boundary condition is set to be zero for the limiting diffusion equation~\eqref{eqn:rho_limit} with the coefficients $b_i$ computed from~\eqref{eqn:f_L_i}.

Similarly we derive the limiting boundary condition for $f^\RR$. Denote
\begin{equation}\label{eqn:f_R_decompose}
f^\RR=\rho(x=1)f_3+\partial_x\rho(x=1)f_4\,,
\end{equation}
with $f_{3,4}$ satisfying reflective boundary layer equation:
\begin{equation}\label{eqn:f_R_i}
\begin{cases}
v\partial_zf_i&=\frac{\langle \mathsf{Kn}\rangle}{\mathsf{Kn}}(\mathcal{L}f_i-f_i)\\
f_i(z=0,v,\omega)&=\eta f_i(z=0,-v,\omega)+\psi_i\quad v>0\\
f_i(z\to\infty)&\in\Null(\mathcal{L}-\mathbb{I})
\end{cases}
\end{equation}
where $\psi_i$, the Dirichlet component of the reflective boundary layer condition, takes the value of
\begin{equation*}\label{eqn:f0}
\psi_3(z=0,v,\omega)=-(1-\eta)\,,\quad \psi_4(z=0,v,\omega)=-(1+\eta)v\Kn\,.
\end{equation*}
Using the same notation as in~\eqref{eqn:b_i}, recalling the zero limiting condition for $f^\RR$, we require
\[
\rho(x=1)b_3+\partial_x\rho(x=1)b_4=0\,.
\]

Summarizing the derivations above, we have the interior solution being~\eqref{eqn:f_interior} with $\rho$ satisfying:
\begin{equation}\label{eqn:f_interior_diff}
\begin{aligned}
    \partial_{xx}\rho&=0\\
    b_1\rho-b_2\partial_x\rho&=b_0, \quad\text{at}\quad x=0\,,\\
    b_3\rho+b_4\partial_x\rho&=0, \quad\text{at}\quad x=1\,,
\end{aligned}
\end{equation}
where $b_i$ are constants defined in~\eqref{eqn:b_i} with $f_i$ solving~\eqref{eqn:f_L_i} and~\eqref{eqn:f_R_i}.
This is a diffusion equation with Robin type boundary condition and is uniquely solvable.

\section{Computation of the diffusion limit}\label{sec:Diffusion limit}
In this section we propose our algorithm that solves the limiting diffusion equation~\eqref{eqn:f_interior_diff} with a fine-tuned boundary condition to respect second order accuracy. The computation of the equation itself is fairly straightforward, and can be done by a standard finite difference method. The crucial part is to find the constants $b_i$. This amounts to finding a proper solver for~\eqref{eqn:f_L_i} and~\eqref{eqn:f_R_i}, as summarized in Algorithm~\ref{alg:linear}.

\RestyleAlgo{boxruled}
\begin{algorithm}
\SetAlgoLined
\KwData{ Kinetic boundary conditions given by \begin{itemize}\item[1.]  $\phi$ for $v>0$ at $x=0$  ;\item[2.] reflection coefficient $\eta$ for $v<0$ at $x=1$\end{itemize}}
\KwResult{Temperature profile from diffusion approximation}
Step I: Compute boundary data\\
\For{$i=0,1,2$}{
at $x=0$: Compute $b_i$ from solving equation~\ref{eqn:f_L_i}.
}
\For{$i=3,4$}{
at $x=1$: Compute $b_i$ from solving equation~\ref{eqn:f_R_i}.
}
Step II: Compute $\rho$ using~\ref{eqn:f_interior_diff}.
\caption{Compute the asymptotic diffusion limit.}\label{alg:linear}
\end{algorithm}

We now separately discuss the computation of~\eqref{eqn:f_L_i},~\eqref{eqn:f_R_i} and~\eqref{eqn:f_interior_diff}. In~\cite{li2017convergent}, the authors proposed a spectral method that computes the solution to these half space boundary layer equation with Dirichlet boundary condition quickly, and it was later extended in~\cite{li2017half} to deal with general boundary data, including reflective or diffusive. Equation~\eqref{eqn:f_L_i} and~\eqref{eqn:f_R_i} are equipped with Dirichlet and reflective boundary conditions respectively, and the solvers can be easily modified to be used here.

\noindent\textbf{Computation of~\eqref{eqn:f_L_i}:} According to~\eqref{eqn:b_i}, $b_i$ is the limit of~\eqref{eqn:f_L_i} with different incoming data $\psi_i$. For the conciseness, we omit the subindex, and solve the following half-space equation with incoming data:
\begin{equation}\label{eqn:f_half_space}
\begin{cases}
v\partial_xf&=\frac{\langle \mathsf{Kn}\rangle}{\mathsf{Kn}}(\mathcal{L}f-f)\\
f(x=0,v,\omega)&=\psi, \quad v>0\\
f(x\to\infty)&=\theta_\infty\in\Null(\mathcal{L}-\mathbb{I})\,.
\end{cases}
\end{equation}
There are two main challenges: the equation experiences a jump discontinuity at $(x=0,v=0)$, and the equation is supported on the full domain, with one end unknown ($\theta_\infty$). To overcome the first difficulty, we will build the jump discontinuity directly into the basis function, by separating $f$ into its even and odd components. The odd part of $f$ naturally encodes the jump discontinuity. To overcome the second difficulty, we will employ a spectral method and decompose $f$ into summation of separable functions. The function on the velocity domain form an orthonormal basis for the solution to expand on. This translates a PDE into a coupled ODE system for coefficients on $x$, and this ODE-system can be solved semi-analytically. Since the solution at $x=\infty$ is not known, we utilize the damping-recovering process, as suggested by Proposition 3.6 in~\cite{li2017half}, to force the final data to be zero. This is to solve the damped equation twice, once with the original boundary condition $\psi$, and another time with $1\in\Null({\mathcal{L}-\mathbb{I}})$, both with zero data at $x=\infty$. To be more specific, we let:
\begin{equation}\label{eqn:f_damped}
\begin{cases}
v\partial_x\Tilde{f}&=\frac{\langle \mathsf{Kn}\rangle}{\mathsf{Kn}}\mathcal{L}_d\Tilde{f}\\
\Tilde{f}(x=0,v,\omega)&=\psi\quad v>0\\
\Tilde{f}(x\to\infty)&=0
\end{cases}
\end{equation}
where the damped operator $\mathcal{L}_d$ is given by
\begin{equation}\label{eqn:damped_operator}
    \mathcal{L}_d\Tilde{f}=(\mathcal{L}-\mathbb{I})\Tilde{f}+\alpha v\langle v,\Tilde{f}\rangle+\alpha v(\mathcal{L}-\mathbb{I})^{-1}v\langle v(\mathcal{L}-\mathbb{I})^{-1}v,\Tilde{f}\rangle
\end{equation}
where $0<\alpha\ll 1$ is the damping parameter. $g_0$ solves the same equation~\eqref{eqn:f_damped} with $\psi$ replaced by $1$ for all positive velocity. According to Proposition 3.6 in~\cite{li2017half}, the solution to~\eqref{eqn:f_half_space} should be:
\begin{equation}\label{eqn:theta_infty_recovery}
f = \tilde{f}-\theta_\infty (g_0-1)\,,\quad\text{where}\quad
   \theta_\infty=\frac{\langle v,\Tilde{f}(0,v,\omega)\rangle}{\langle v,g_0(0,v,\omega)\rangle}\,.
\end{equation}

Having translated the computation of~\eqref{eqn:f_half_space} to the computation of the damped equation~\eqref{eqn:f_damped}, we now solve this damped equation using even-odd spectral decomposition. More specifically, let $\tilde{f}^E$ and $\tilde{f}^O$ denote the even part and odd part of $\tilde{f}$ respectively on the velocity space:
\begin{equation}
    \tilde{f}^E(x,v,\omega)=\frac{\tilde{f}(x,v,\omega)+\tilde{f}(x,-v,\omega)}{2},\quad\text{and}\quad \tilde{f}^O(x,v,\omega)=\frac{\tilde{f}(x,v,\omega)-\tilde{f}(x,-v,\omega)}{2}\,.
\end{equation}
Then immediately we have:
\[
\tilde{f} = \tilde{f}^E+\tilde{f}^O\,,\quad\text{with}\quad \tilde{f}^E(v) = \tilde{f}^E(-v)\,,\quad\text{and}\quad \tilde{f}^O(v) = -\tilde{f}^O(-v)\,.
\]
For the even and odd parts of the solution, we use the Legendre polynomial basis functions $\{\l_k\}_{k=1}^N$ to expand them out such that
\begin{equation*}
    \int\limits_{0}^1\l_i(v) \l_j(v) \mathrm{d}v=\delta_{ij}\,.
\end{equation*}
In order to account for $\omega$, we give a frequency dependent weight to these Legendre polynomials such that
\begin{equation*}
    \phi_i(v,\omega)=\sqrt{\frac{C_\tau}{C_\omega/\tau}}l_i(v)\,.
\end{equation*}
This definition then is combined with the even or odd extension to the full velocity domain:
\[
\phi^E_i(v,\omega) = \begin{cases} \phi_i(v,\omega)\,,\quad v>0\\ \phi_i(-v,\omega)\, \quad v<0\end{cases}\,,\quad\text{and}\quad \phi^O_i(v,\omega) = \begin{cases} \phi_i(v,\omega)\,,\quad v>0\\ -\phi_i(-v,\omega)\, \quad v<0\end{cases}\,.
\]

Note that this expansion defines a set of orthonormal basis:
\[
\frac{1}{2}\langle \phi^E_m \phi^O_n\rangle=0\,,\quad  \frac{1}{2}\langle \phi^E_m \phi^E_n\rangle=\delta_{mn}\,,\quad\text{and}\quad \frac{1}{2}\langle \phi^O_m \phi^O_n\rangle=\delta_{mn}\,,
\]
where the bracket notation means
\begin{equation}\label{eqn:bracket}
    \langle f\rangle=\int\limits_{v=-1}^1\int\limits_{\omega\ge 0}\frac{C_\omega/\tau}{C_\tau}f\mathrm{d}\omega\mathrm{d}v\,.
\end{equation}
It is on this basis we solve the damped equation~\eqref{eqn:f_damped}:
\begin{equation}\label{eqn:f_tilde_c}
\tilde{f} = \tilde{f}^E+\tilde{f}^O = \sum_{m}c^E_m\phi^E_m+ \sum_{m}c^O_m\phi^O_m\,.
\end{equation}

Substituting this into the damped equation gives 
\begin{equation}\label{eqn:generalized_e_value1}
    A\frac{\rd}{\rd x} \Vec{c}=B\Vec{c}\,,
\end{equation}
where $\Vec{c}$ contains the weights $\{c_k^E\}_{k=1}^N, \{c_k^O\}_{k=1}^{N+1}$ and $A, B$ are given by
\begin{equation}
    A=
    \begin{pmatrix}
    A^v & 0\\
    0 & A^v
    \end{pmatrix},\quad
    B=
    \begin{pmatrix}
    0 & B^E\\
    B^O & 0
    \end{pmatrix},\quad\text{and}\quad
    c=
    \begin{pmatrix}
    c^O\\
    c^E
    \end{pmatrix}\,.
\end{equation}
where 
\begin{equation}
    A^v_{ij}=2\langle v\phi_i(v,\omega)\phi_j(v,\omega)\rangle\,,\quad B^E_{ij}=\frac{\langle\Kn\rangle}{\Kn}\langle\phi^E_i,\mathcal{L}_d\phi^E_j\rangle,\quad\text{and}\quad B^O_{ij}=\frac{\langle\Kn\rangle}{\Kn}\langle\phi^O_i,\mathcal{L}_d\phi^O_j\rangle\,.
\end{equation}

The equation~\eqref{eqn:generalized_e_value1} is not yet solvable. To uniquely determine this, we utilize two sources of information.
\begin{itemize}
    \item There are $N$ constraints imposed on the incoming boundary. This is to set, for all $j=1,\cdots,N$:
    \begin{equation}\label{eqn:N_constraints1}
\begin{aligned}
    &\sum\limits_{i=1}^{N+1}c^O_i\langle v\phi_j^E\phi^O_i\rangle_++\sum\limits_{i=1}^{N}c^E_i\langle v\phi_j^E\phi^E_i\rangle_+=\langle v\phi_j^E \psi\rangle_+\,.
\end{aligned}
\end{equation}
Here $\langle\cdot\rangle_+$ is similarly defined as in~\eqref{eqn:bracket} except the velocity domain is confined for $v>0$.
    \item The other $N+1$ constraints come from the well-posedness requirement. To ensure that the solution $\Vec{c}(x)\to0$ as $x\to\infty$, $\Vec{c}$ projected on the ``growing" modes of~\eqref{eqn:generalized_e_value1} should be all eliminated at the incoming boundary location. To identify the ``growing" mode, we consider the following generalized eigenvalue problem.
\begin{equation}\label{eqn:generalized_eigenvalue}
    \lambda_kA\Vec{v}_k=B\Vec{v}_k\,.
\end{equation}
Using the transformation $\Vec{e}_k=\Vec{v}_k^\top A \Vec{c}$, we have the following equation for $\Vec{e}_k$.
\begin{equation}
    \partial_x\Vec{e}_k=\lambda_k\Vec{e}_k\,.
\end{equation}
This suggests model $\Vec{e}_k$ will either exponentially increase or decrease as $x\to\infty$. To ensure the solution vanishes at infinity, we tolerate arbitrary projection of the solution on modes $\Vec{e}_k$ that decrease, but require the projection on positive and zero modes ($\lambda_k\geq 0$) to be zero at $x=0$. According to~\cite{li2017convergent}, there are exactly $N+1$ such modes, and we have the remaining $N+1$ constraints in the form of
\begin{equation}\label{eqn:c_0_decay_cond}
  \Vec{e}_k(0)=  \Vec{v}_k^\top A\Vec{c}(0)=0\,,\quad\text{for}\quad \lambda_k\geq 0\,.
\end{equation}
\end{itemize}
The two equations,~\eqref{eqn:N_constraints1} and~\eqref{eqn:c_0_decay_cond} combined, uniquely determine $\Vec{c}(x=0)$, making~\eqref{eqn:generalized_e_value1} solvable with diminishing solution at $x=\infty$. $\tilde{f}$ is also determined by~\eqref{eqn:f_tilde_c}.

One solves~\eqref{eqn:f_damped} again using $\psi$ replaced by $1$ for $g_0$, and compute $\theta_\infty$ using~\eqref{eqn:theta_infty_recovery}.

\noindent\textbf{Computation of~\eqref{eqn:f_R_i}:} Similarly, to compute $b_{3,4}$, we need to compute the half-space equation with reflective boundary condition~\eqref{eqn:f_R_i} using different $\psi_i$. We supress the subscript, and compute the following:
\begin{equation}\label{eqn:f_damped_reflective}
\begin{cases}
v\partial_xf&=\frac{\langle \mathsf{Kn}\rangle}{\mathsf{Kn}}(\mathcal{L}f-f)\\
f(x=0,v,\omega)&=\eta f(x=0,-v,\omega)+\psi\quad v>0\\
f(x\to\infty)&=\theta_\infty\in\Null(\mathcal{L}-\mathbb{I})
\end{cases}\,.
\end{equation}
We now follow the damping-recovering process in Proposition 3.6 in~\cite{li2017convergent} for:
\[
f=\tilde{f}-\theta_\infty (g_0-1)\,,
\]
where $\tilde{f}$ solves the following damped equation:
\begin{equation}
\begin{cases}
v\partial_x\Tilde{f}&=\frac{\langle \mathsf{Kn}\rangle}{\mathsf{Kn}}\mathcal{L}_d\Tilde{f}\\
\Tilde{f}(x=0,v,\omega)&=\eta\Tilde{f}(x=0,-v,\omega)+\psi\quad v>0\\
\Tilde{f}(x\to\infty)&=0
\end{cases}\,,
\end{equation}
and $g_0$ solves the same damped equation above with the source term given by $1-\eta$. Here,~\ref{eqn:N_constraints1} is modified as, for all $j=1\,,\cdots,N$: 
\begin{equation}
\begin{aligned}
    &\sum\limits_{i=1}^{N+1}c^O_i\langle v\phi_j^E\phi^O_i\rangle_++\sum\limits_{i=1}^{N}c^E_i\langle v\phi_j^E\phi^E_i\rangle_+ =\eta\Big[-\sum\limits_{i=1}^{N+1}c^O_i\langle v\phi_j^E\phi^O_i\rangle_+ +\sum\limits_{i=1}^{N}c^E_i\langle v\phi_j^E\phi^E_i\rangle_+ \Big]+\langle v\phi_j^E \psi\rangle_+\,.
\end{aligned}
\end{equation}
The projection on the ``growing"-mode is set to be zero, as is done for the incoming boundary condition. Repeat the process for computing $g$, one finally gets:
\begin{equation}\label{eqn:theta_infty_r}
   \theta_\infty=\langle v,\Tilde{f}(0,v,\omega)\rangle/\langle v,g_0(0,v,\omega)\rangle\,.
\end{equation}

\noindent\textbf{Computation of~\eqref{eqn:f_interior_diff}:} To solve $\rho$, we use the standard finite difference scheme. Let the domain $[0,1]$ be discretized using uniform grid points $\{x_i\}_{i=0}^{N_x}$ with grid spacing $\Delta x=\frac{1}{N_x}$, then numerically, we solve for $\rho(x_i)=\rho_i, i=0,\cdots N_x$ as
\begin{equation}\label{eqn:rho_discrete}
\begin{cases}
        \frac{1}{(\Delta x)^2}(\rho_{i+1}-2\rho_{i}+\rho_{i-1})&=0, \quad i=1,\cdots, N_x-1\\
        b_1\rho_0-\frac{b_2}{\Delta x}(\rho_1-\rho_0)&=b_0\\
        b_3\rho_{N_x}+\frac{b_4}{\Delta x}(\rho_{N_x}-\rho_{N_x-1})&=0
\end{cases}
\end{equation}
where $b_i, i=0,\cdots,4$ are prepared from above. We note that with ghost-cell strategy, higher order discretization is also possible. This is no longer the main point of the paper and is thus omitted.
\section{Numerical Results}\label{sec:Numerical results}
We showcase two numerical examples to demonstrate the asymptotic convergence in this section.

\begin{figure}[!h]
 \begin{subfigure}{0.5\textwidth}     
    \centering
    \includegraphics[width=7cm]{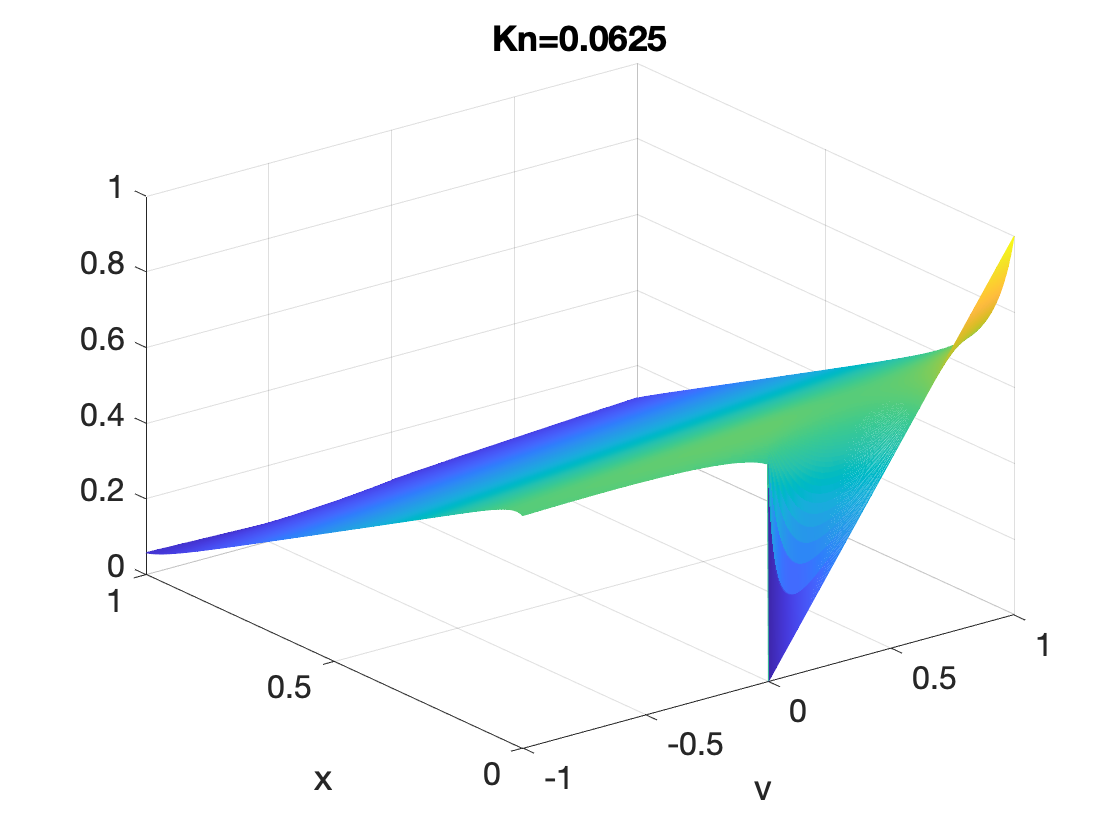}
    \caption{$\phi=v$}
   \end{subfigure}
 \begin{subfigure}{0.5\textwidth}     
    \centering
    \includegraphics[width=7cm]{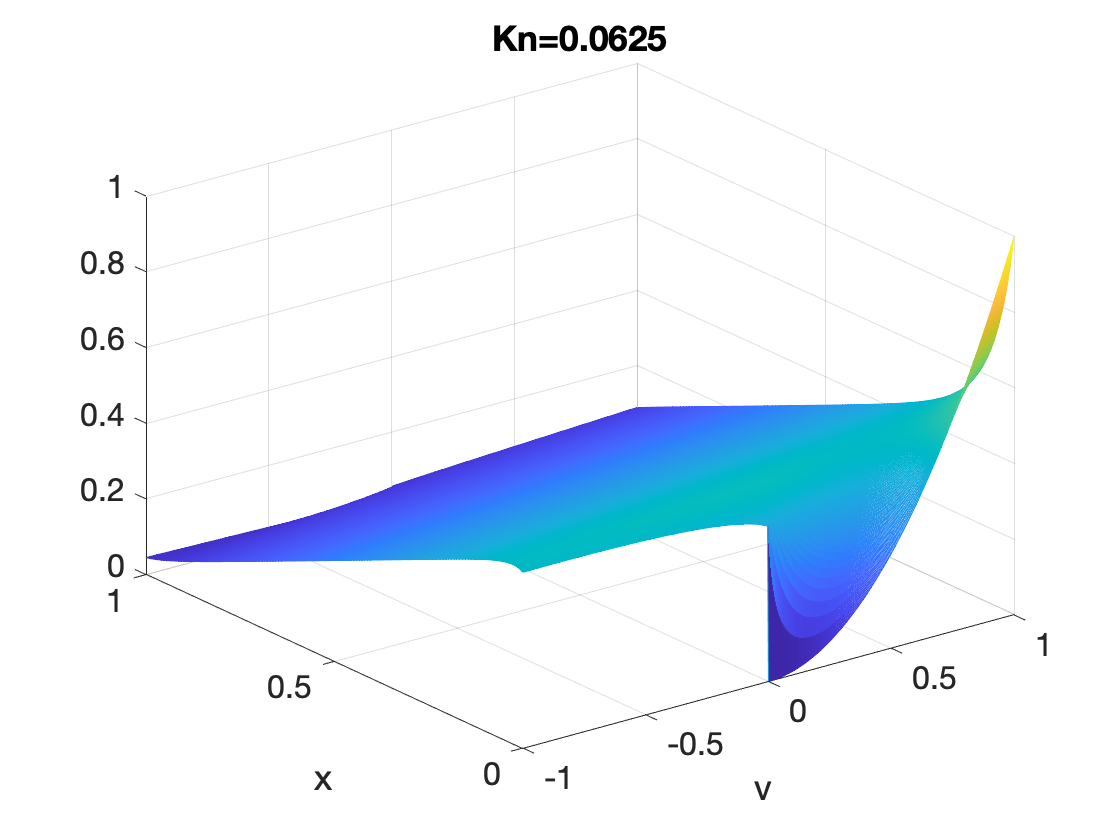}
    \caption{$\phi=v^2$}
    \end{subfigure}
        \caption{Example I. Reference solution for $\Kn=0.0625$.}
        \label{fig:reference_f}
    \end{figure}

\noindent\textbf{Example I.} We first demonstrate the numerical result when there is only one frequency in the system, namely we assume $f(x,v,\omega)=f(x,v)\delta_{\omega-1}$. We let $\eta=0.5$ in this case, and compute the reference solution of~\eqref{eqn:f_original} using $\Delta x = \Delta v = 2^{-9}$ and $\alpha = 0.01$. This reference solution is computed using the standard finite volume with upwinding fluxes. The smallest $\Kn$ used in our computation is $\Kn = 1/16$ so the discretization is fine enough to resolve the layer. We show the reference solution using $\phi=v$ and $\phi=v^2$ for $\Kn = 1/16$ as the incoming data in Figure~\ref{fig:reference_f}.

\begin{figure}[!h]
    \centering
    \includegraphics[width = 0.3\textwidth,height= 0.15\textheight]{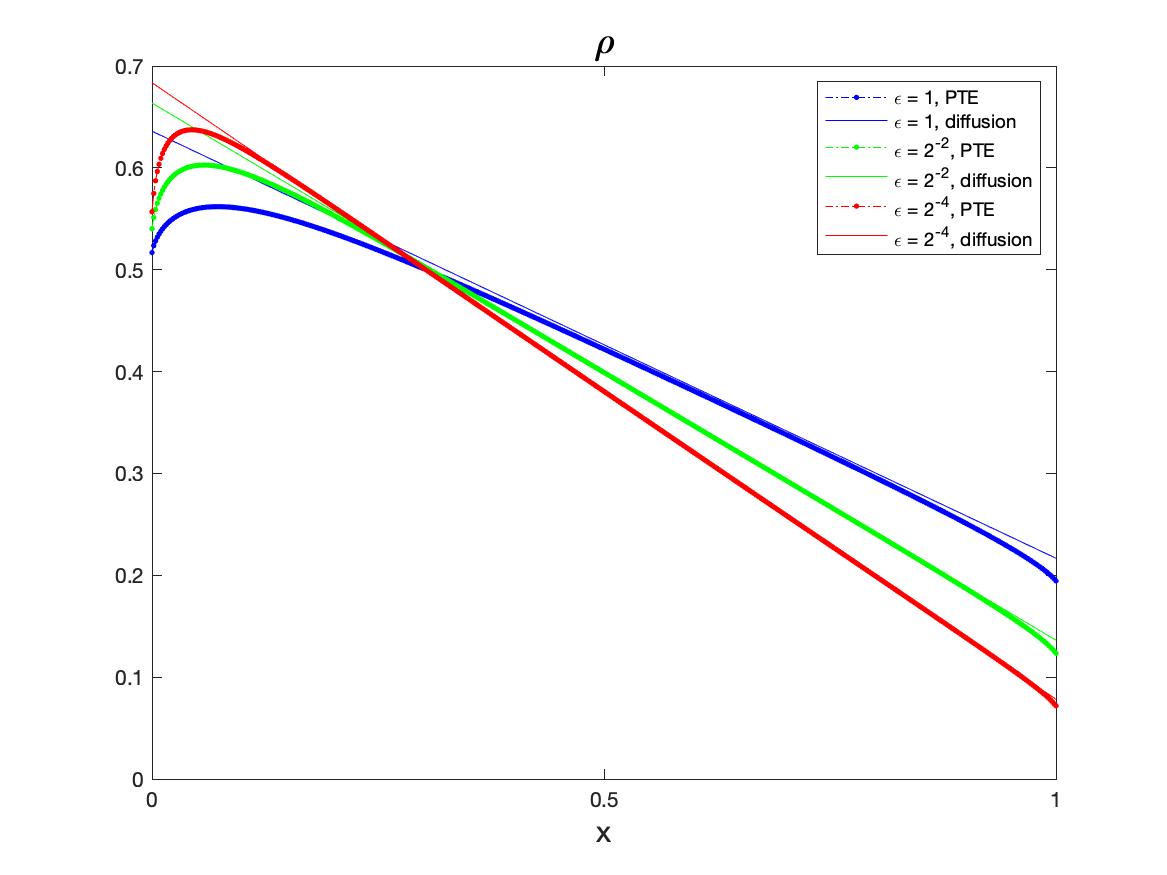}
    \includegraphics[width = 0.3\textwidth,height= 0.15\textheight]{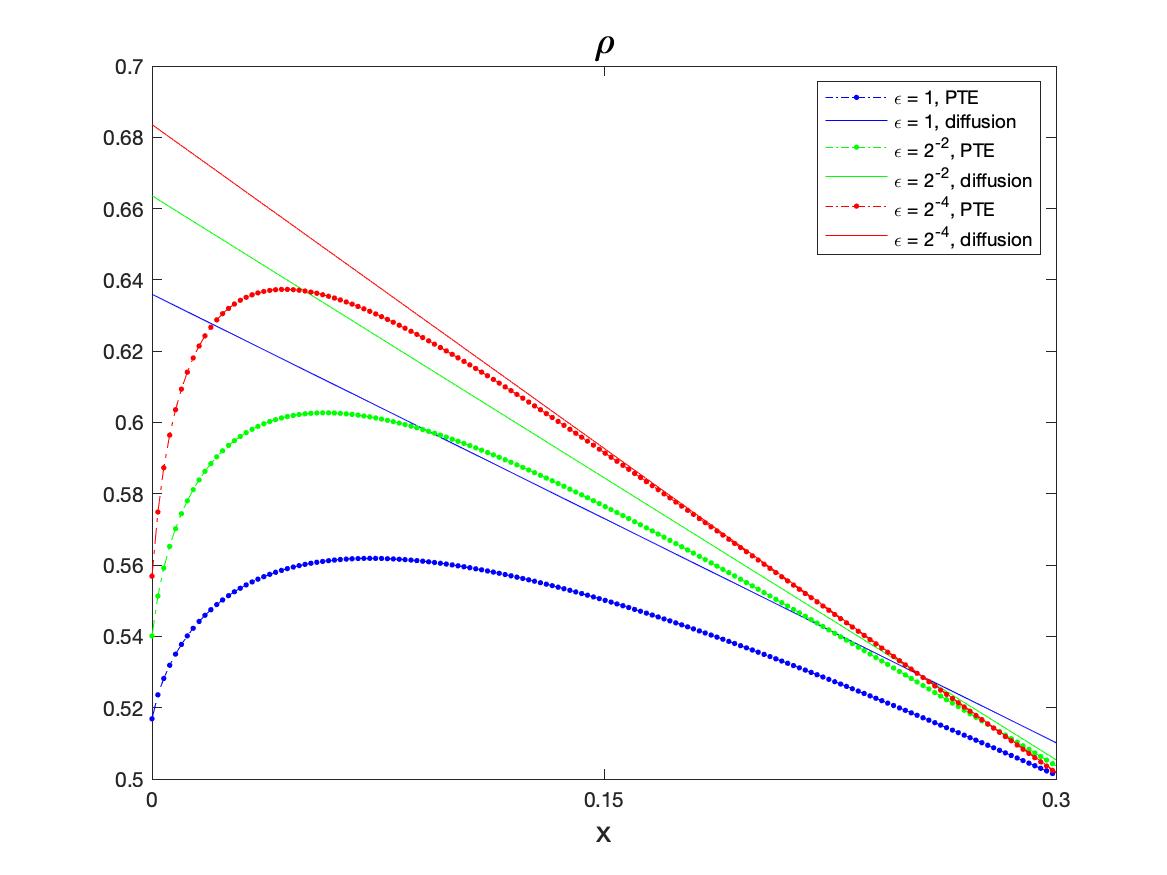}
    \includegraphics[width = 0.3\textwidth,height= 0.15\textheight]{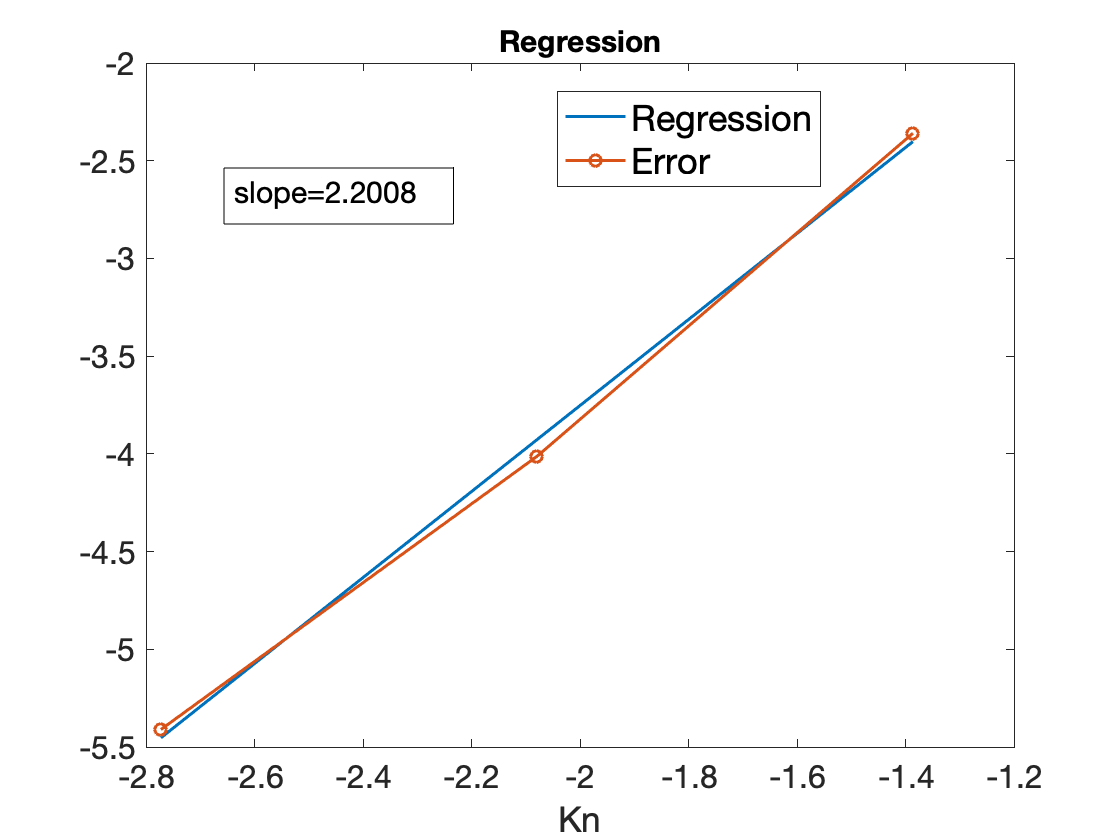}
    \caption{Example I. The panel on the left shows the density $\rho$ over the whole domain. The panel in the middle shows the layer behavior close to $x=0$ computed using different $\Kn$ and the limiting $\rho$. The panel on the right shows the convergence rate on the log-log scale. It suggests the asymptotic convergence is $\Kn^2$. The incoming data is $\phi = v$.}
    \label{fig:robin_both_sides_one_frequency_ex1}
\end{figure}

\begin{figure}[!h]
\centering
    \includegraphics[width = 0.3\textwidth,height= 0.15\textheight]{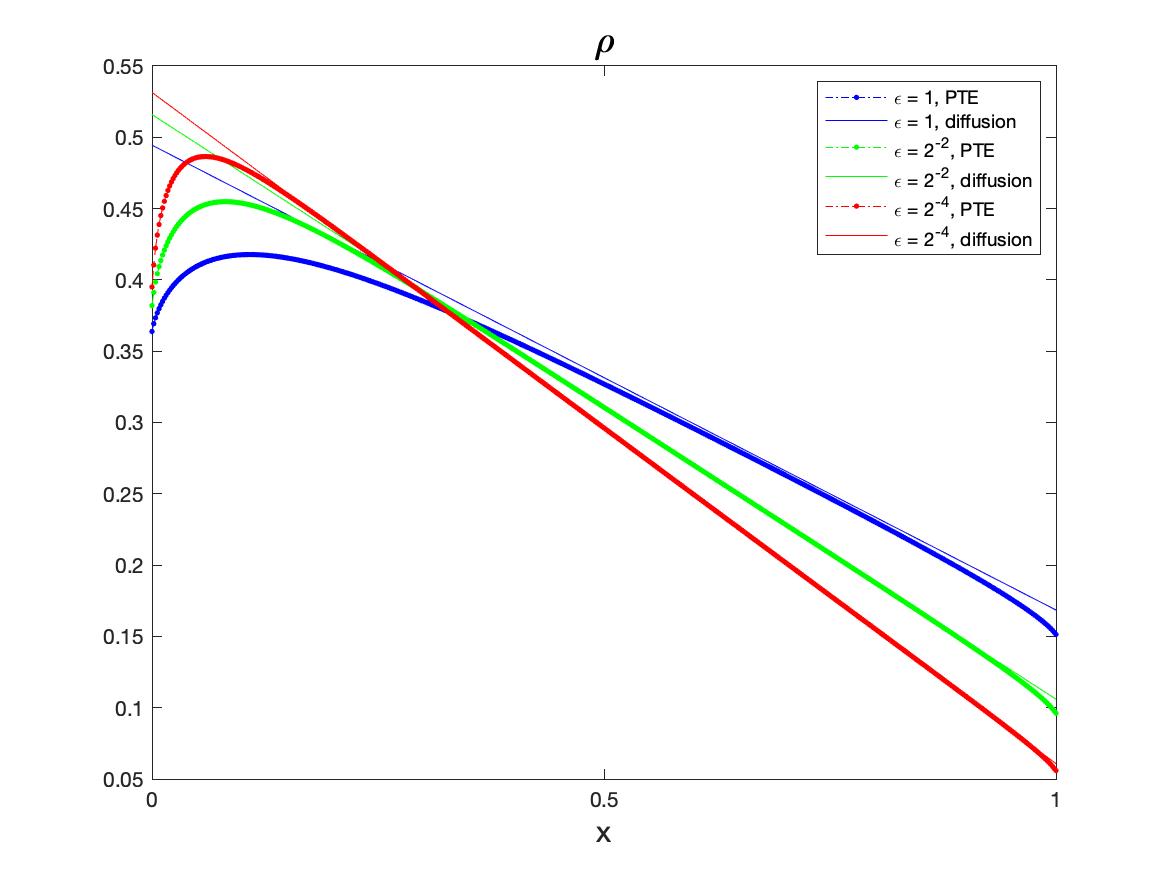}
    \includegraphics[width = 0.3\textwidth,height= 0.15\textheight]{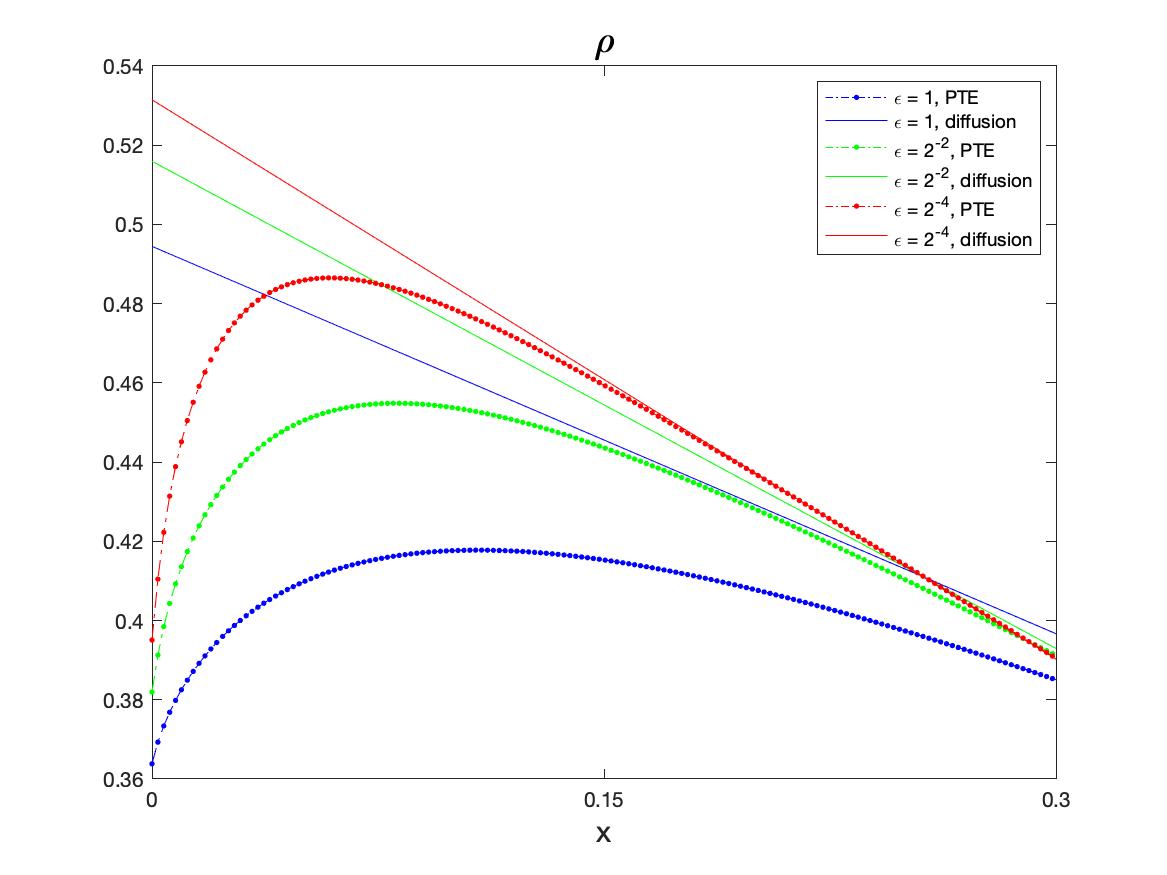}
    \includegraphics[width = 0.3\textwidth,height= 0.15\textheight]{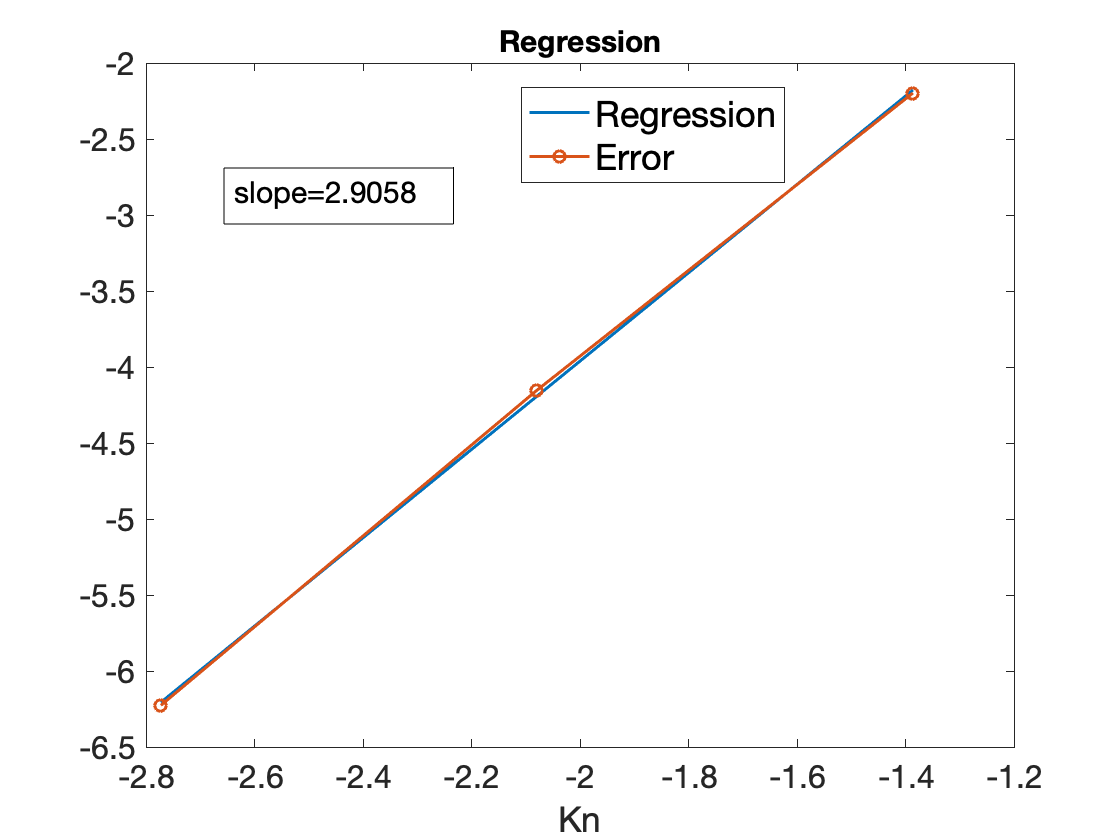}
\caption{Example I. The panel on the left shows the density $\rho$ over the whole domain. The panel in the middle shows the layer behavior close to $x=0$ computed using different $\Kn$ and the limiting $\rho$. The panel on the right shows the convergence rate on the log-log scale. It suggests the asymptotic convergence is $\Kn^2$. The incoming data is $\phi = v^2$.}
    \label{fig:robin_both_sides_one_frequency_ex2}
\end{figure}

The asymptotic approximation solution is computed using~\eqref{eqn:f_interior_diff} also using $\Delta x=2^{-9}$. In Figure~\ref{fig:robin_both_sides_one_frequency_ex1}, we plot the comparison between the kinetic solution and its diffusion limit for different $\Kn$. As $\Kn$ changes from $1/4$ to $1/8$ and finally to $1/16$, the layer effect is more and more obvious. In the first panel we show the convergence rate in $\Kn$, and it is clear the convergence rate is $2$, indicating the asymptotic error is $\Kn^2$. Here the error is defined as:
\[
\text{Error}=\sqrt{\sum_{i=[N_x/4]}^{[3N_x/4]}|\rho_i-T_i|^2}
\]
where $\rho_i$ is computed from the diffusion equation with Robin boundary condition, and $T_i$ denotes the value of temperature at $x_i$ computed using the reference kinetic equation. Clearly, the error excludes the layer fluctuation. The same plots are generated for $\phi = v^2$ as well, see Figure~\ref{fig:robin_both_sides_one_frequency_ex2}. From the plots, the layer effect is obvious, and as $\Kn$ decreases, the layers become more and more sharp, and the kinetic solutions with boundary layers get closer and closer to the limiting diffusion equation equipped with Robin boundary condition. We should note that the Robin boundary condition itself has $\Kn$ dependence, so the profile for the diffusion limit differ for different $\Kn$.

\begin{figure}[!h]
 \begin{subfigure}{0.5\textwidth}
    \centering
    \includegraphics[width=7cm]{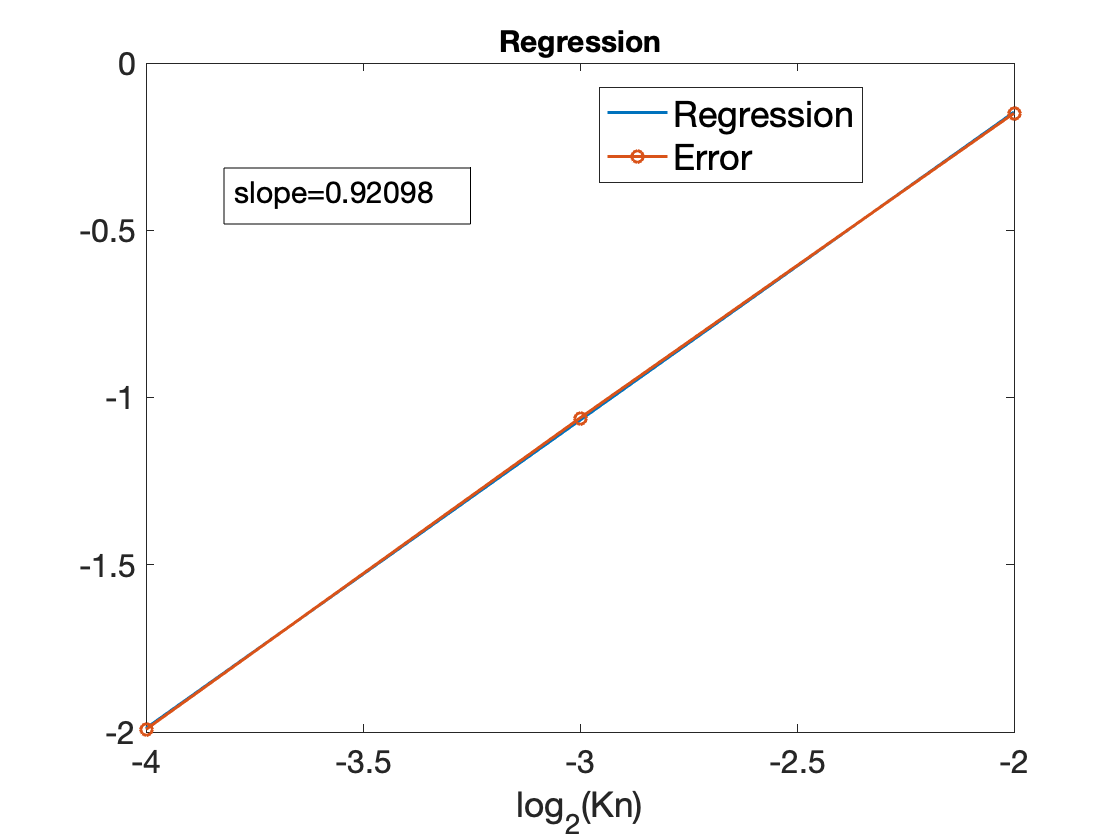}
    \caption{$\phi=v$}
 \end{subfigure}
 \begin{subfigure}{0.5\textwidth}
    \centering
    \includegraphics[width=7cm]{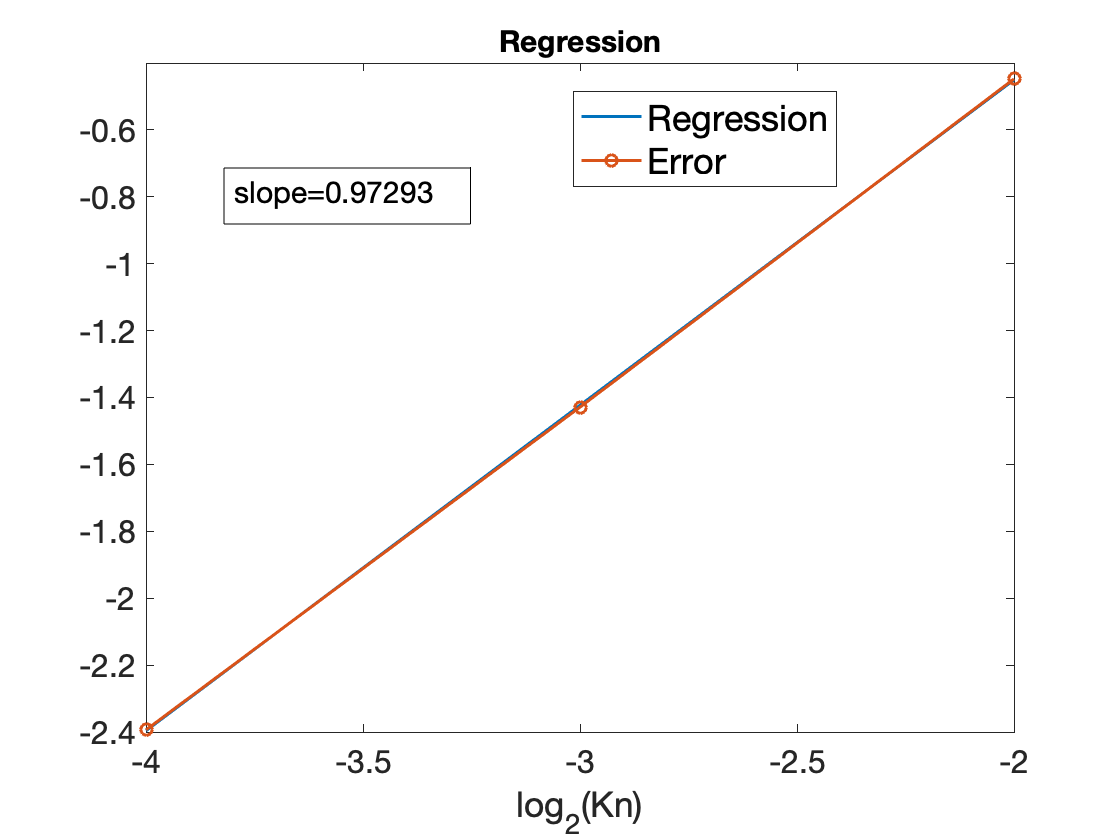}
    \caption{$\phi=v^2$}
        \end{subfigure}
            \caption{Example I. If Dirichlet boundary condition is used for the limiting equation, first order convergence is obtained.}
            \label{fig:O(Kn)_one_freq}
    \end{figure}

We should mention that our result also recovers that in~\cite{li2015diffusion}. If Dirichlet-type boundary condition is used for computing $\rho$, instead of Robin-type, then the asymptotic convergence rate will degrade to the first order. This is reflected in Figure~\ref{fig:O(Kn)_one_freq}. For both $\phi=v$ and $\phi=v^2$, we use the Dirichlet boundary condition at $x=0$ to compute the limit and show the error convergence with respect to $\Kn$ on the log-log scale. The convergence rate is approximately $\sim 1$, as expected.

\begin{figure}[!h]
 \begin{subfigure}{0.5\textwidth}     
    \centering
    \includegraphics[width=7cm]{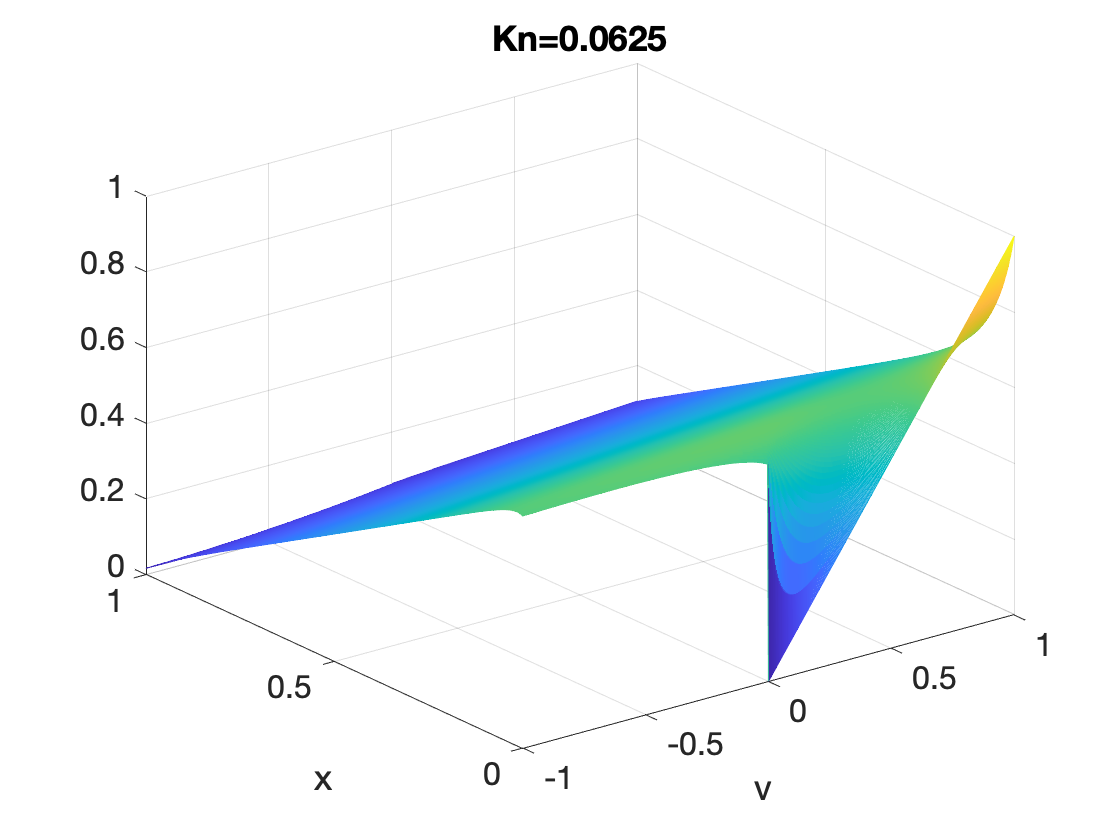}
        \caption{$f(x,v,\omega)$ for $\omega=1.2$}
   \end{subfigure}
 \begin{subfigure}{0.5\textwidth}     
    \centering
    \includegraphics[width=7cm]{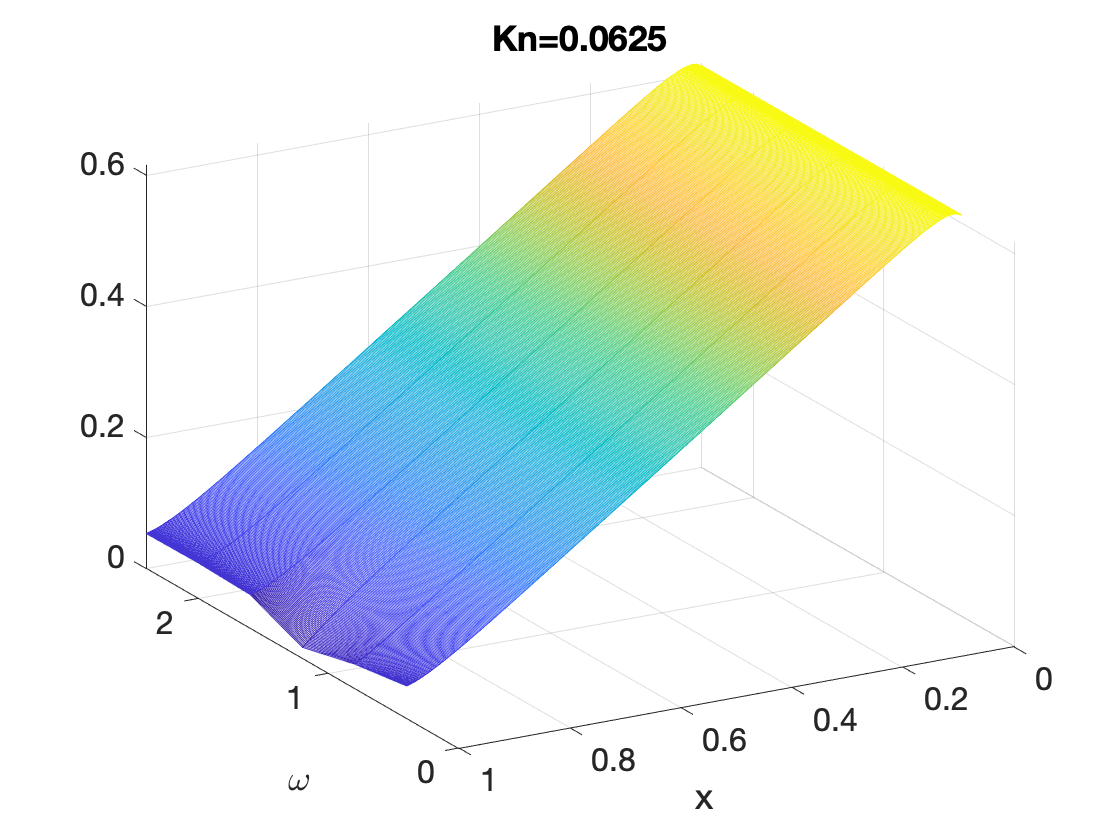}
            \caption{$f(x,v,\omega)$ for $v=-1$}
    \end{subfigure}
        \caption{Reference solution for $\Kn=0.0625$ in the multi-frequency case when $\phi=v$.}
        \label{fig:reference_f_multi_freq_phi_v}
    \end{figure}

\begin{figure}[!h]
    \centering
    \includegraphics[width = 0.3\textwidth,height= 0.15\textheight]{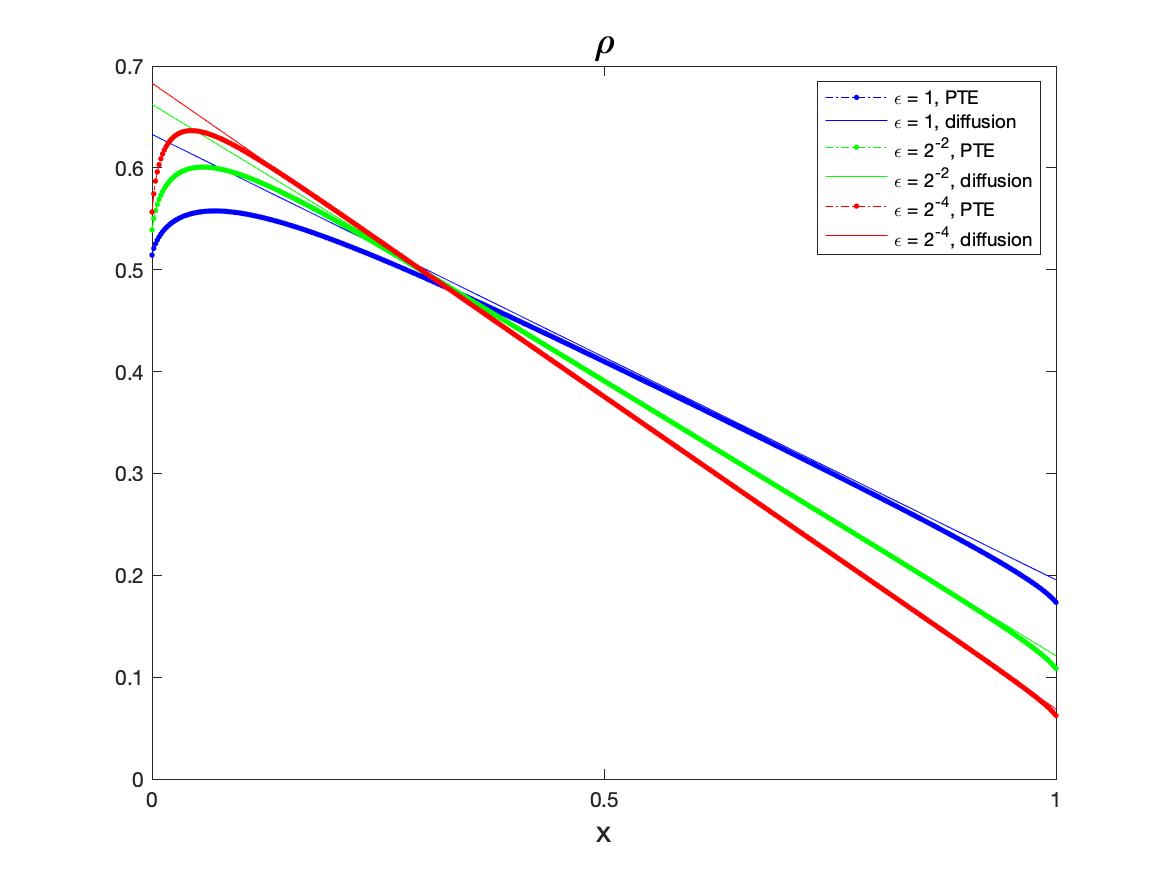}
    \includegraphics[width = 0.3\textwidth,height= 0.15\textheight]{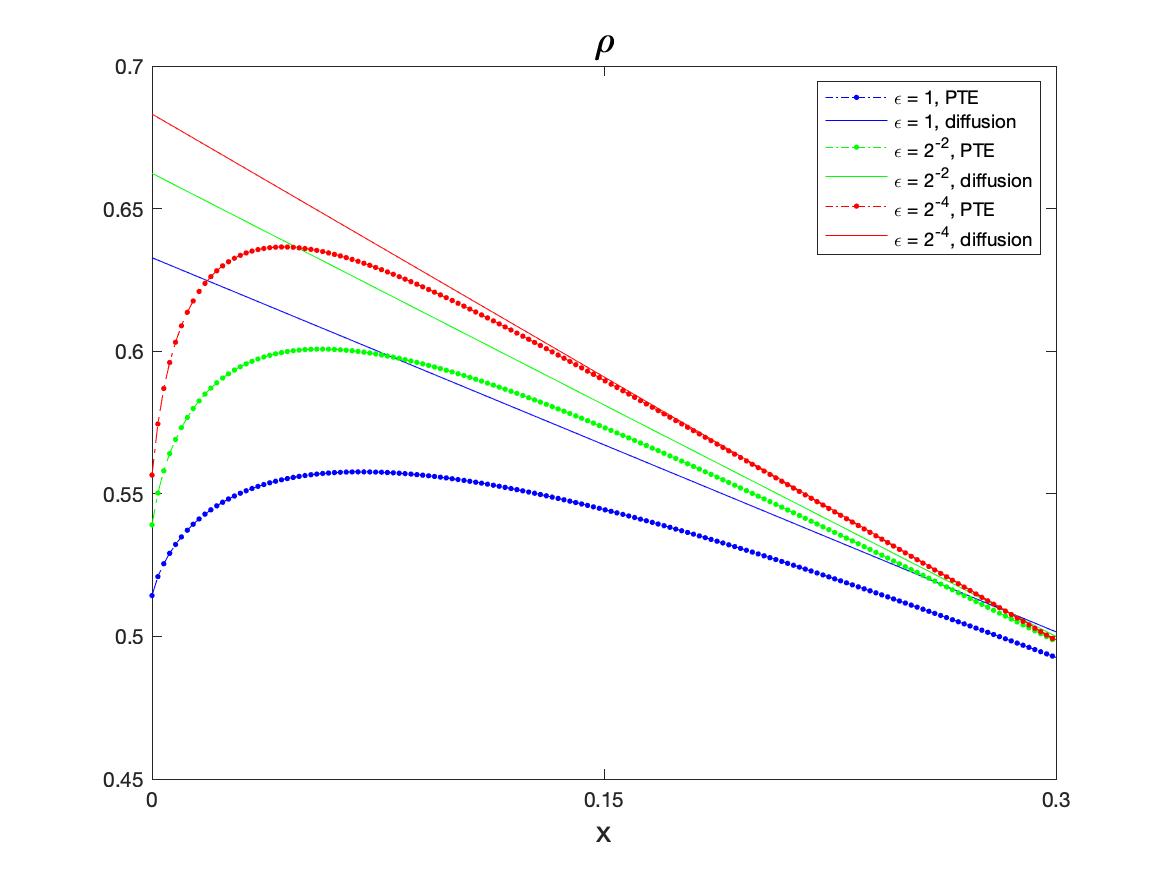}
    \includegraphics[width = 0.3\textwidth,height= 0.15\textheight]{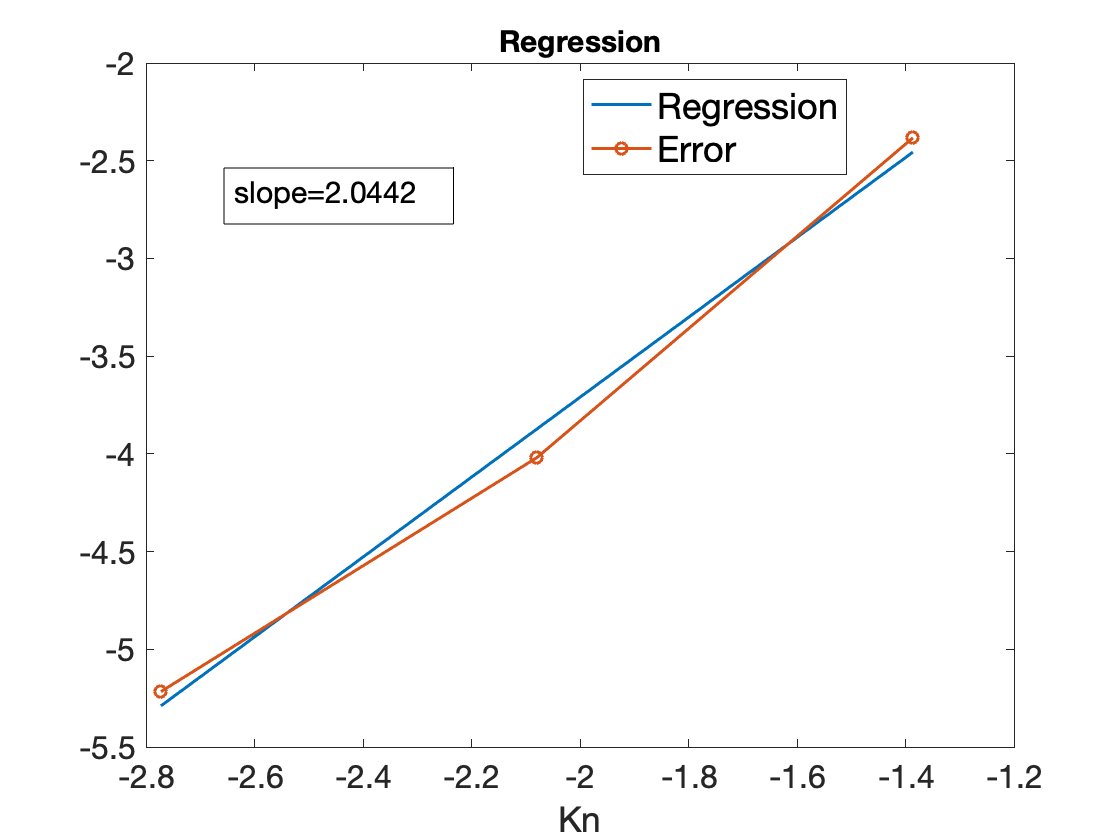}
\caption{Example II. The panel on the left shows the density $\rho$ over the whole domain. The panel in the middle shows the layer behavior close to $x=0$ computed using different $\Kn$ and the limiting $\rho$. The panel on the right shows the convergence rate on the log-log scale. It suggests the asymptotic convergence is $\Kn^2$. The incoming data is $\phi = v$ for the multiple frequency.}
\label{fig:robin_both_sides_multi_frequency_ex1}
\end{figure}

\begin{figure}[!h]
    \centering
    \includegraphics[width = 0.3\textwidth,height= 0.15\textheight]{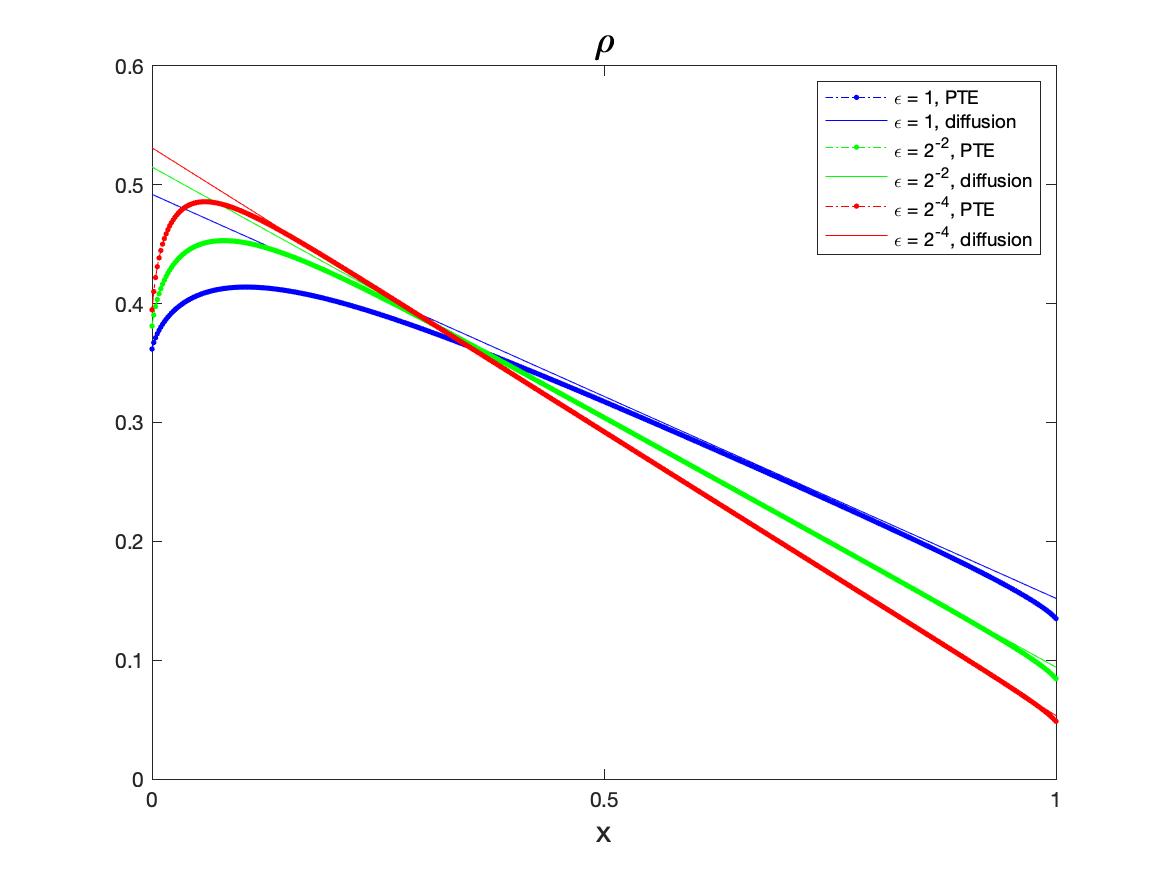}
    \includegraphics[width = 0.3\textwidth,height= 0.15\textheight]{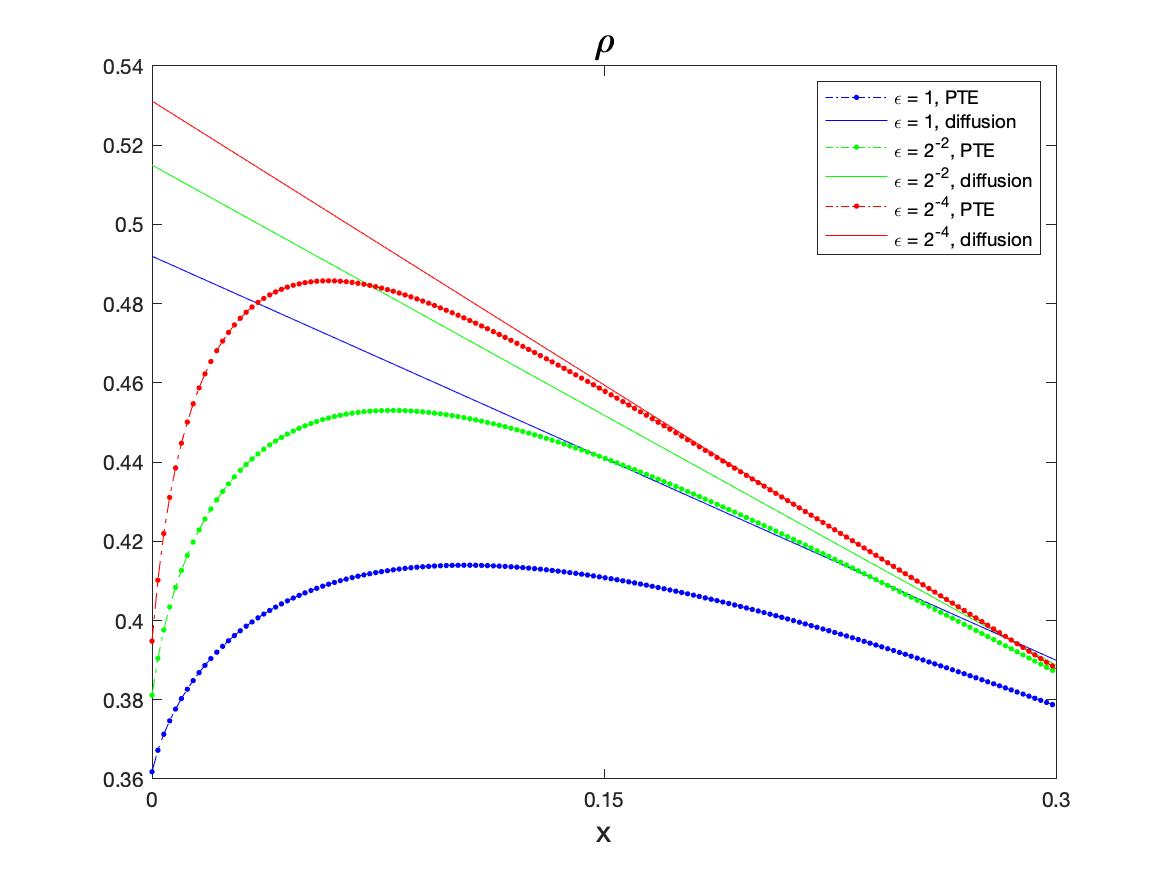}
    \includegraphics[width = 0.3\textwidth,height= 0.15\textheight]{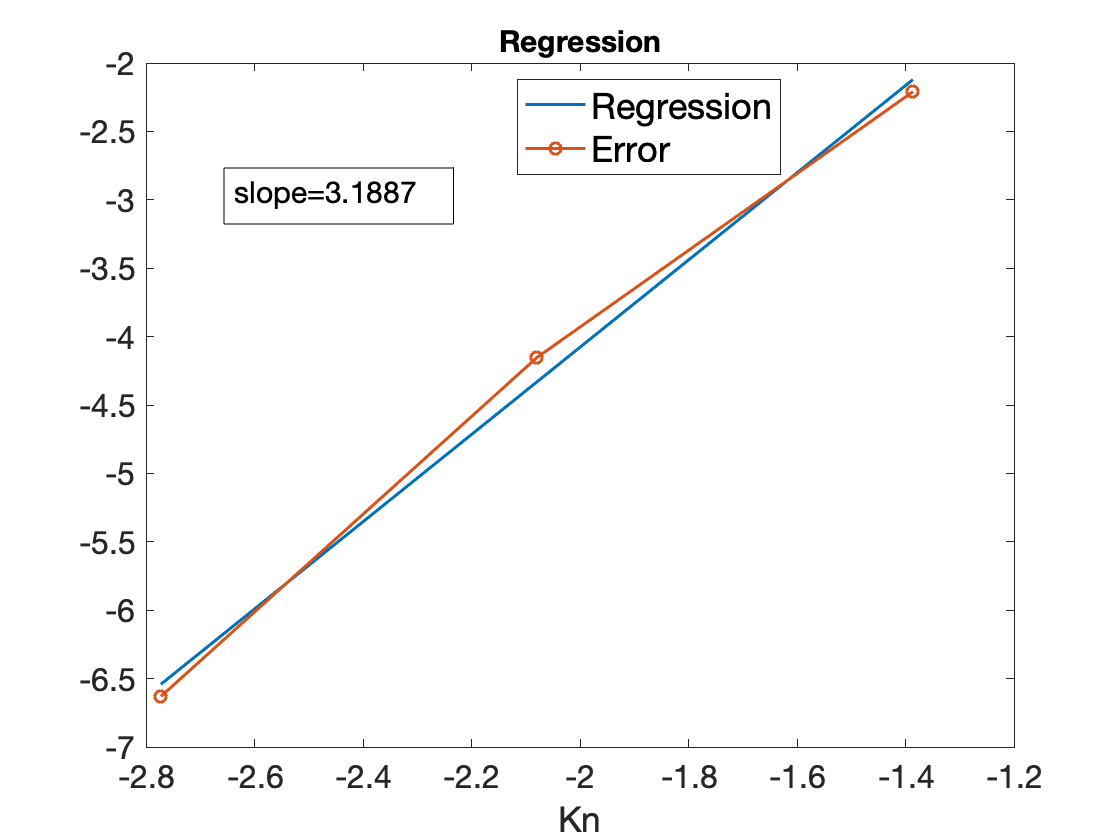}
\caption{Example II. The panel on the left shows the density $\rho$ over the whole domain. The panel in the middle shows the layer behavior close to $x=0$ computed using different $\Kn$ and the limiting $\rho$. The panel on the right shows the convergence rate on the log-log scale. It suggests the asymptotic convergence is $\Kn^2$. The incoming data is $\phi = v^2$ for the multiple frequency.}
\label{fig:robin_both_sides_multi_frequency_ex2}
\end{figure}

\begin{figure}[!h]
 \begin{subfigure}{0.5\textwidth}
    \centering
    \includegraphics[width=7cm]{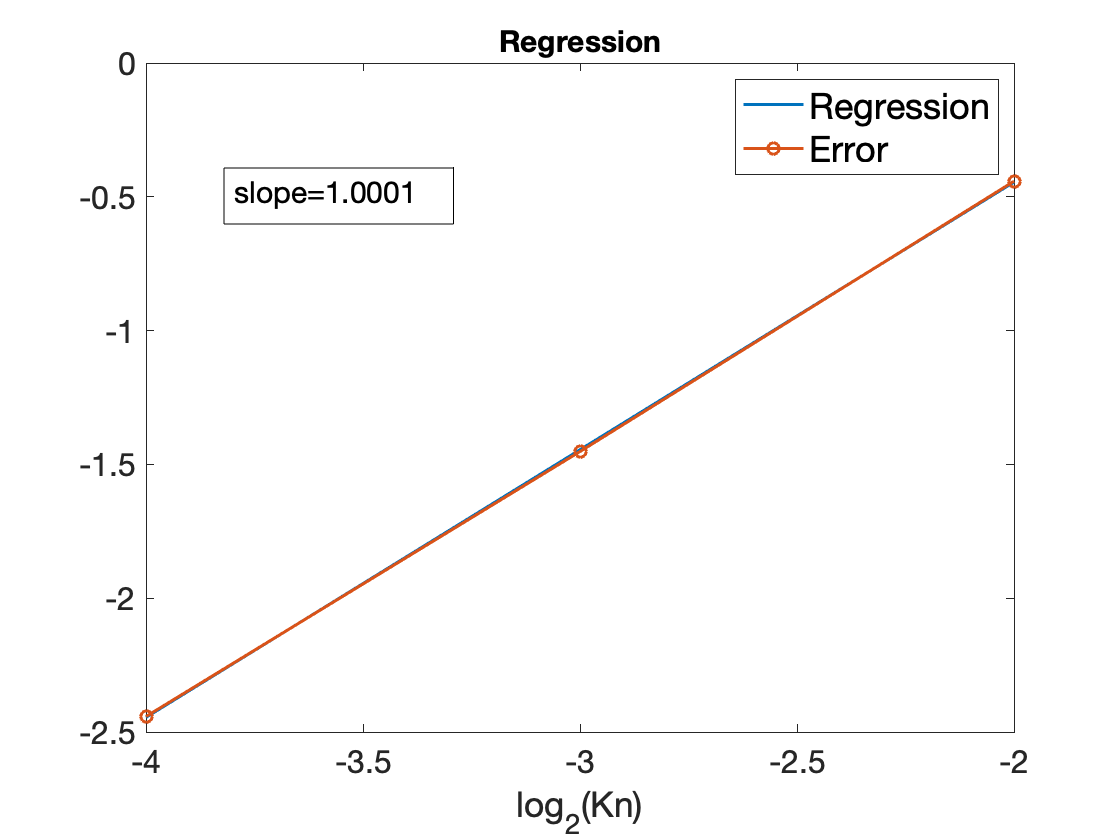}
    \caption{$\phi=v$}
 \end{subfigure}
 \begin{subfigure}{0.5\textwidth}
    \centering
    \includegraphics[width=7cm]{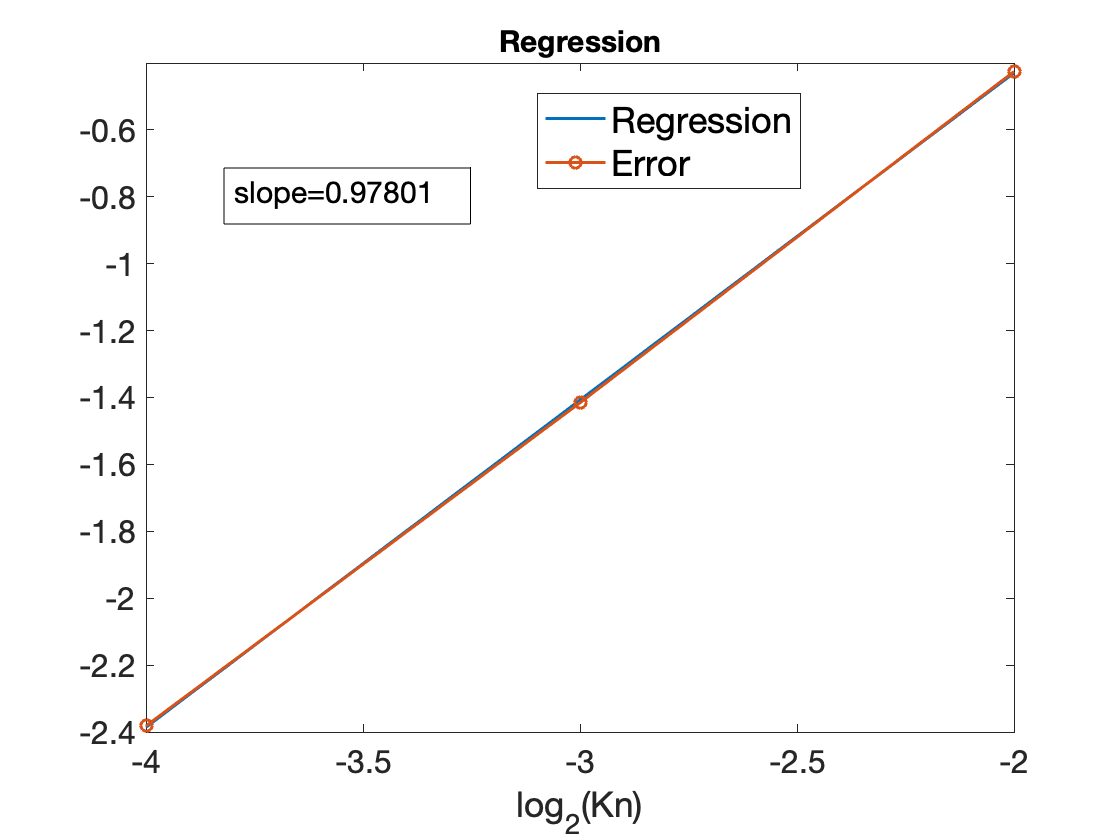}
    \caption{$\phi=v^2$}
        \end{subfigure}
            \caption{We show the comparison between the error vs $\mathrm{Kn}$ plot in the interior and the best fit using linear regression using the $\mathcal{O}(\Kn)$ scheme for the multi-frequency case.}
            \label{fig:O(Kn)_multi_freq}
    \end{figure}

\noindent\textbf{Example II.} The second example concerns the case where there are multiple frequencies. We uniformly discretize the frequency domain $\omega\in[0.4,2.4]$ into six bins with grid spacing being $\Delta\omega = 0.4$. The reflection coefficient is defined using the following:
\begin{equation}
    \eta(\omega)=)= \frac{1}{2}+\frac{\tanh(10(\omega-1.5))-\tanh(2(\omega-1))}{4}\,.
\end{equation}
Taking $D_\omega =1$, we have, in this case, $C_\omega$, as defined in~\eqref{eqn:C_omega} taking the form of $C_\omega=\frac{(10\omega)^2 e^{10\omega}}{\big(e^{10\omega}-1\big)^2}$. We also choose $\tau(\omega)=\frac{1}{10\omega}$ and $\|\boldsymbol v_g\|=10\omega$ as suggested in~\cite{forghani2018reconstruction}. The reference solution is once again computed using $\Delta x = \Delta v = 2^{-9}$ and $\alpha = 0.01$, using the standard finite volume with upwinding fluxes. One reference solution computed using $\Kn = 1/16$ is shown in Figure~\ref{fig:reference_f_multi_freq_phi_v}. The asymptotic limiting equation is computed according to~\eqref{eqn:f_interior_diff} with Robin boundary condition. Both for $\phi = v$ and $\phi=v^2$, the layer become sharper as $\Kn\to0$, and the convergence is $O(\Kn^2)$, as shown in Figure~\ref{fig:robin_both_sides_multi_frequency_ex1} and~\ref{fig:robin_both_sides_multi_frequency_ex2}. As in the single frequency case, the layer effect is obvious, and small $\Kn$ brings the kinetic solution closer to the diffusion solution. If Dirichlet boundary is used in place of Robin, then the convergence rate deteriorates to the first order. This is shown in Figure~\ref{fig:O(Kn)_multi_freq}.

\noindent\textbf{Acknowledgement:} The research of Q.L. and A.N. is supported in part by NSF via grant DMS-1750488 and Office of the Vice Chancellor for Research and Graduate Education at the University of Wisconsin Madison with funding from the Wisconsin Alumni Research Foundation. W.S. acknowledges support from NSERC Discovery Grant R611626. 

\noindent\textbf{Data availability statement:} The authors are happy to provide programming files and data if requested.

\newpage
\bibliographystyle{unsrt}
\bibliography{Reference}
\end{document}